\documentclass[10pt,a4paper]{article}
\usepackage[]{geometry}

\title{A Roadmap for the Computation of Persistent Homology}

\usepackage{array} 
\usepackage[caption]{subfig} 
\usepackage{hyperref,enumitem}
\usepackage{amsmath}
\usepackage{rotating}
\usepackage{amsfonts}
\usepackage{multicol} 
\usepackage{framed}
\usepackage{algorithm,algorithmic}
\usepackage{cite}
\usepackage{graphicx,rotating,booktabs,framed,placeins}
\usepackage{epsfig}

\usepackage{enumitem}
\usepackage{verbatim}
\usepackage{cprotect}
 \usepackage{alltt}
 \usepackage{breakurl}
 \usepackage{tablefootnote}

\usepackage{adjustbox}
\usepackage{multirow}

\usepackage{standalone,tikz}
\usetikzlibrary{patterns}
\newcommand*\MyScale{1}

\tikzset{every picture/.style={scale=\MyScale}}

\newtheorem{thm}{Theorem}
\newtheorem{definition}[thm]{Definition}

\newtheorem{rmk}[thm]{Remark}

\newcommand{\javaplex}{\textsc{javaPlex}}
\newcommand{\perseus}{\textsc{Perseus}}
\newcommand{\gudhi}{\textsc{Gudhi}}

\newcommand{\jholes}{\textsc{jHoles}}
\newcommand{\dionysus}{\textsc{Dionysus}}
\newcommand{\phat}{\textsc{PHAT}}
\newcommand{\dipha}{\textsc{DIPHA}}
\newcommand{\plex}{\textsc{Plex}}
\newcommand{\jplex}{\textsc{jPlex}}
\newcommand{\simppers}{\textsc{SimpPers}}
\newcommand{\tdapackage}{\textsc{TDA Package}}
\newcommand{\rivet}{\textsc{RIVET}}
\newcommand{\ripser}{\textsc{ripser}}
\newcommand{\hera}{\textsc{hera}}
\newcommand{\landscape}{\textsc{Persistence Landscape Toolbox}}

\newcommand{\klein}{\mathrm{\textbf{Klein}}}

\newcommand{\random}{\mathrm{\textbf{random}}}

\newcommand{\genome}{\mathrm{\textbf{genome}}}
\newcommand{\eleg}{\mathrm{{\it \textbf{eleg}}}}
\newcommand{\HIV}{\mathrm{\textbf{HIV}}}
\newcommand{\influenza}{\mathrm{\textbf{H3N2}}}
\newcommand{\dragone}{\mathrm{\textbf{drag 1}}}
\newcommand{\dragtwo}{\mathrm{\textbf{drag 2}}}
\newcommand{\fract}{\mathrm{\textbf{fract}}}
\newcommand{\fractr}{\mathrm{\textbf{fract r}}}

\newcommand{\netwsc}{\mathrm{\textbf{netw-sc}}}
\newcommand{\vicsek }{\mathrm{\textbf{Vicsek}}}
\newcommand{\vertebra}{\mathrm{\textbf{vertebra}}}
\newcommand{\senate}{\mathrm{\textbf{senate}}}
\newcommand{\house}{\mathrm{\textbf{house}}}

\newcommand{\cpp}{\textsf{C++}}

\newcommand{\java}{\textsf{Java}}
\newcommand{\python}{\textsf{Python}}

\newcommand{\gcc}{\textsf{gcc}}
\newcommand{\openmpi}{\textsf{openmpi}}

\newcommand{\VR}{\mathrm{VR}}
\newcommand{\cech}{\check{C}}
\newcommand{\del}{\mathrm{Del}}
\newcommand{\W}{\mathrm{W}}
\newcommand{\WRCF}{\mathrm{WRCF}}

\newcommand{\transpose}{\mathrm{t}}

\newcommand{\low}{\mathrm{low}}

\newcommand{\image}{\mathrm{image}}
\newcommand{\kernel}{\mathrm{kernel}}
\newcommand{\bigo}{\mathcal{O}}

\newcommand{\collection}{\mathcal{U}}
\newcommand{\dist}{\mathrm{d}}
\newcommand{\dg}{\mathrm{dg}}

\usepackage{wrapfig}
\graphicspath{{figures/}} 
\usepackage{booktabs} 
\usepackage[font=small,labelfont=bf]{caption} 
\usepackage{tikz}
\usetikzlibrary{shadows,arrows,positioning}
\usepackage{booktabs} 
\pgfdeclarelayer{background}
\pgfdeclarelayer{foreground}
\pgfsetlayers{background,main,foreground}
 \usepackage{textcomp}
 \usepackage{pgfplots}
\tikzstyle{materia}=[draw, fill=blue!20, text width=6.0em, text centered,
  minimum height=1.5em]
  \tikzstyle{material}=[draw,  text width=40.0em,
  minimum height=4em]
    \tikzstyle{mater}=[draw, fill=blue!20, text width=3.8em, text centered,
  minimum height=3em]
\tikzstyle{etape} = [materia, text width=4.5em, minimum width=2em,
  minimum height=3em, rounded corners]
  \tikzstyle{etap} = [mater, text width=5em, minimum width=3em,
  minimum height=3em, rounded corners]
 \tikzstyle{etapee} = [materia, text width=6.2em, minimum width=5.5em,
  minimum height=3em, rounded corners]

\tikzstyle{texto} = [above, text width=8em, text centered]
\tikzstyle{linepart} = [draw, thick, color=black!50, -latex', dashed]
\tikzstyle{line} = [draw,-, thick, color=black!50, -latex', dashed]
\tikzstyle{arrow} = [draw,->, ultra thick, color=black!500, -latex']
\tikzstyle{ur}=[draw, text centered, minimum height=0.01em]

\newcommand{\etape}[2]{node (p#1) [etape]
  {#2}}
  
\newcommand{\etapee}[2]{node (p#1) [etapee]
  {#2}}

\renewcommand{\thefootnote}{\arabic{footnote}}
\author{Nina Otter ({\tt otter@maths.ox.ac.uk})\footnotemark[1]\ \footnotemark[3]
\and Mason~A. Porter ({\tt mason@math.ucla.edu})\footnotemark[1]\ \footnotemark[3]\ \footnotemark[4]\ 
\and Ulrike Tillmann ({\tt tillmann@maths.ox.ac.uk})\footnotemark[1]\ \footnotemark[3]
\and Peter Grindrod ({\tt grindrod@maths.ox.ac.uk})\footnotemark[1]\
\and Heather~A. Harrington ({\tt harrington@maths.ox.ac.uk})\footnotemark[1]\
}

\date{}
\begin{document}
\maketitle

\renewcommand{\thefootnote}{\fnsymbol{footnote}}

\footnotetext[1]{Mathematical Institute, University of Oxford, Oxford OX2 6GG, UK}
\footnotetext[2]{CABDyN Complexity Centre, University of Oxford, Oxford OX1 1HP, UK}
\footnotetext[3]{The Alan Turing Institute, 96 Euston Road, London NW1 2DB, UK}
\footnotetext[4]{Department of Mathematics, UCLA, Los Angeles, CA 90095, USA}

\renewcommand{\thefootnote}{\arabic{footnote}}


 \begin{abstract}
 
 \emph{Persistent homology} (PH) is a method used in topological data analysis  (TDA) to study qualitative features of data that persist across multiple scales. It is robust to perturbations of input data, independent of dimensions and coordinates, and provides a compact representation of the qualitative features of the input. 
The computation of PH is an open area with numerous important and fascinating challenges. 
 The field of PH computation is evolving rapidly, and new algorithms and software implementations are being updated and released at a rapid pace. The purposes of our article are to (1)  introduce theory and computational methods for PH to a broad range of 
computational scientists and (2) provide benchmarks of state-of-the-art implementations for the computation of PH.  
We give a friendly introduction to PH, navigate the  for the computation of PH with an eye towards applications, and use a range of synthetic and real-world data sets to evaluate currently available open-source implementations for the computation of PH. Based on our benchmarking, we indicate which algorithms and implementations are best suited to different types of data sets. In an accompanying  tutorial, we provide guidelines for the computation of PH. We make publicly available all scripts that we wrote for the tutorial, and we make available the processed version of the data sets used in the benchmarking.

 \end{abstract}

\pagestyle{myheadings}
\thispagestyle{plain}


\subsection*{Keywords: persistent homology, topological data analysis, point-cloud data, networks}

\newpage

\section{Introduction} \label{intro}

The amount of available data has increased dramatically in recent years, and this situation ---  which will only become more extreme ---  necessitates the development of innovative and efficient data-processing methods.
Making sense of the vast amount of data is difficult: on one hand, the sheer size of the data poses challenges;  on the other hand, the complexity of the data, which includes situations in which data is noisy, high-dimensional, and/or incomplete, is perhaps an even more significant challenge. The use of clustering techniques and other ideas from areas such as computer science, machine learning, and uncertainty quantification --- along with mathematical and statistical models --- are often very useful for data analysis (see, e.g., \cite{kaufman1990,airoldi2010,siam-cluster,satu07} and many other references). However, recent mathematical developments are shedding new light on such ``traditional'' ideas,  forging new approaches of their own, and helping people to better decipher increasingly complicated  structure in data.

Techniques from the relatively new subject of ``topological data analysis'' (TDA) have provided a wealth of new insights in the study of data in an increasingly diverse set of applications --- including 
sensor-network coverage \cite{dSG07}, 
proteins \cite{KNBNH14,GHIKMN15,XW14,XLM16}, 
three-dimensional structure of DNA \cite{ESR15},
development of cells \cite{RCKR17},
stability of fullerene molecules \cite{XFTW15}, 
robotics \cite{BGK15,PHR15,VAB13}, 
signals in images \cite{CBK09,GBKP13}, 
periodicity in time series \cite{PDHH15}, 
cancer \cite{NLC11,DCCVPA10,CMCMR16,SCMPN14}, 
phylogenetics \cite{CCR13,CLR16,ERCR14},
natural images \cite{CIdSZ08}, 
the spread of contagions \cite{TKHKMPM14,lo16}, 
self-similarity in geometry \cite{macpherson12},
materials science \cite{KGKM13,KGKM14,hiraoka16,lee17},
financial networks \cite{LPRS08,gidea2017}, 
diverse applications in neuroscience \cite{giusti2016,curto2016,dlotko16,kanari16,lord16,miller16,Yoo2016,dabaghian14}, 
classification of weighted networks \cite{sizemore15},
collaboration networks \cite{Pal2017,CH13},
analysis of mobile phone data \cite{Bajardi2015},
collective behavior in biology \cite{TZH14}, time-series output of dynamical systems \cite{MZR15}, and more. 
There are numerous others, and new applications of TDA appear in journals and preprint servers increasingly frequently. There are also interesting computational efforts, such as \cite{WW16}.

TDA is a field that lies at the intersection of data analysis, algebraic topology, computational geometry, computer science, statistics, and other related areas. The main goal of TDA is to use ideas and results from geometry and topology  to develop tools for studying qualitative features of data. To achieve this goal, one needs precise definitions of qualitative features, tools to compute them in practice, and some guarantee about the robustness of those features. One way to address all three points is a method in TDA called {\it persistent homology} (PH). This method is appealing for applications because it is based on algebraic topology, which gives a well-understood theoretical framework to study qualitative features of data with complex structure,  is computable via linear algebra, and is  robust with respect to small perturbations in input data. 
 
 Types of data sets that can be studied with PH include finite metric spaces, digital images, level-sets of real-valued functions, and networks (see Section \ref{S:data}). In the next two paragraphs, we give some motivation for the main ideas of persistent homology by discussing two examples of such data sets.  
 
Finite metric spaces are also called ``point cloud'' data sets in the TDA literature. From a  topological point of view, finite metric spaces do not contain any interesting information. One thus considers a  ``thickening'' of a  point cloud at different ``scales of resolution'' and then analyses the evolution of the resulting shape across the different resolution scales. 
The qualitative features are given by topological invariants, and one can represent the variation of such invariants across the different resolution scales in a compact way to summarize the ``shape'' of the data. 

\begin{center}
\begin{figure}[h!]
\includegraphics[scale=0.32]{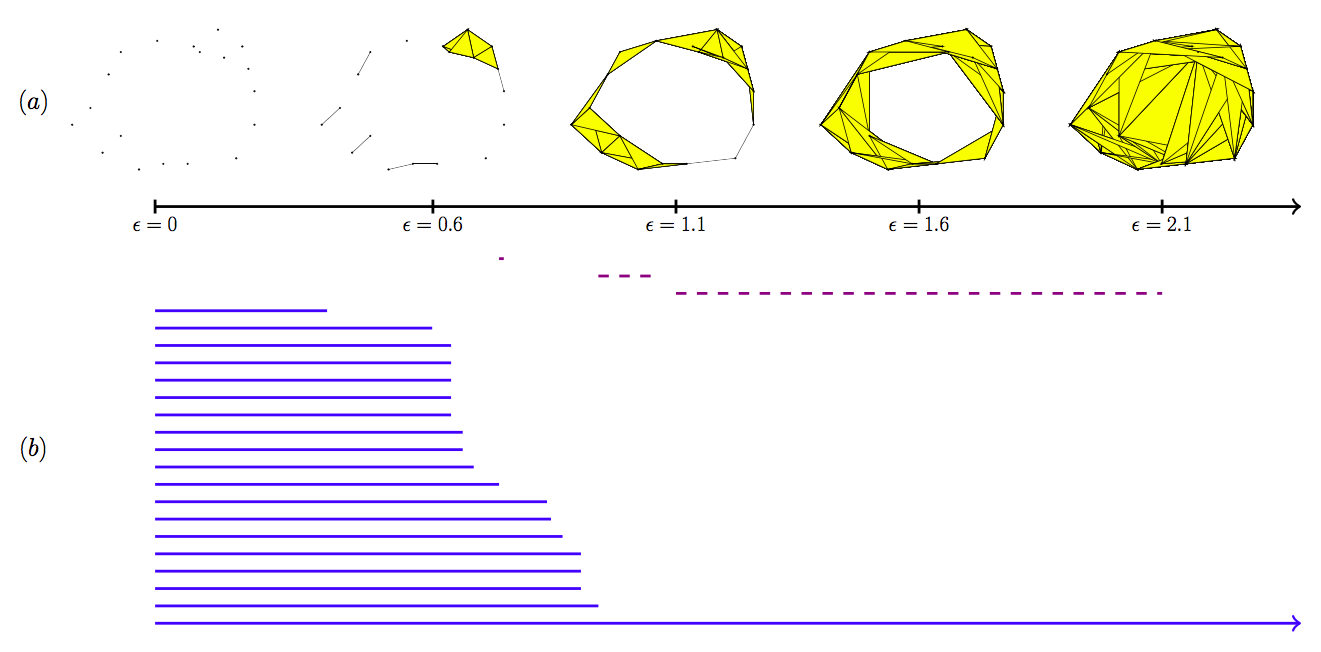}
\caption{
(a) A finite set of points in $\mathbb{R}^2$ (for $\epsilon=0$) and a nested sequence of spaces obtained from it (from $\epsilon=0$ to $\epsilon=2.1$).
(b) Barcode for the nested sequence of spaces illustrated in (a). Solid lines represent the lifetime of components, and dashed lines represent the lifetime of holes.
}\label{figure:filtration}
\end{figure}
\end{center}

 As an illustration, consider the set of points in $\mathbb{R}^2$ that we show in Fig.~\ref{figure:filtration}.
 Let $\epsilon$, which we interpret as a ``distance parameter,'' be a nonnegative real number (so $\epsilon=0$ gives the set of points). For different values of $\epsilon$, we construct a space $S_\epsilon$ composed of vertices, edges, triangles, and higher-dimensional polytopes according to the following rule: We include an edge between two points $i$ and $j$ if and only if the Euclidean distance between them is no larger than $\epsilon$; we include a triangle if and only if all of its edges are in $S_\epsilon$; we include a tetrahedron if and only if all of its face triangles are in $S_\epsilon$; and so on. For $\epsilon\leq \epsilon'$, it then follows that the space $S_\epsilon$ is contained in the space $S_{\epsilon'}$. This yields a nested sequence of spaces, as we illustrate in Fig.~\ref{figure:filtration}(a). Our construction of nested spaces gives an example of a ``filtered Vietoris--Rips complex,'' which we define and discuss in Section \ref{subsection: filtrations}.

By using homology, a tool in algebraic topology, one can measure several features of the spaces $S_\epsilon$ --- including the numbers of components, holes, and voids (higher-dimensional versions of holes). One can then represent the lifetime of such features using a finite collection of intervals known as a ``barcode.'' Roughly, the left endpoint of an interval represents the birth of a feature, and its right endpoint represents the death of the same feature. 
In Fig.~\ref{figure:filtration}(b), we reproduce such intervals for the number of components (blue solid lines) and the number of holes (violet dashed lines).
In Fig.~\ref{figure:filtration}(b), we observe a dashed line that is significantly longer than the other dashed lines. This indicates that the data set has a long-lived hole.  
 By contrast, one can potentially construe the shorter dashed lines as noise. However, note that while widespread, such an interpretation is not correct in general; for applications in which one considers some short and medium-sized intervals as features rather than noise, see ~\cite{SHP16,miller16}. When a feature is still ``alive'' at the largest value of $\epsilon$ that we consider, the lifetime interval is an infinite interval, which we indicate by putting an arrowhead at the right endpoint of the interval. In Fig.~\ref{figure:filtration}(b), we see that there is exactly one solid line that lives up to $\epsilon=2.1$. One can use information about shorter solid lines to extract information about how  data is clustered  in a similar way as with linkage-clustering methods
\cite{siam-cluster}. 

One of the most challenging parts of using PH is statistical interpretation of results.  From a statistical point of view, a barcode like the one in Fig.~\ref{figure:filtration}(b) is an unknown 
quantity that one is trying to estimate; one therefore needs methods for quantitatively assessing the quality of the barcodes that one obtains with computations. The challenge is twofold.  On one hand, there is a cultural obstacle: practitioners of TDA
often have backgrounds in pure topology and are not well-versed in statistical approaches to data analysis \cite{adler14}. On the other hand, the space of barcodes lacks geometric properties that would make it easy to define basic concepts such as mean, median, and so on. Current research is focused both on studying geometric properties of this space and on studying methods that map this space to spaces that have better geometric properties for statistics. 
 In Section \ref{SS:interpretation}, we give a brief overview of the challenges and current approaches for statistical interpretation of barcodes. This is an active area of research and an important endeavor, as
  few statistical tools are currently available for interpreting results in applications of PH.

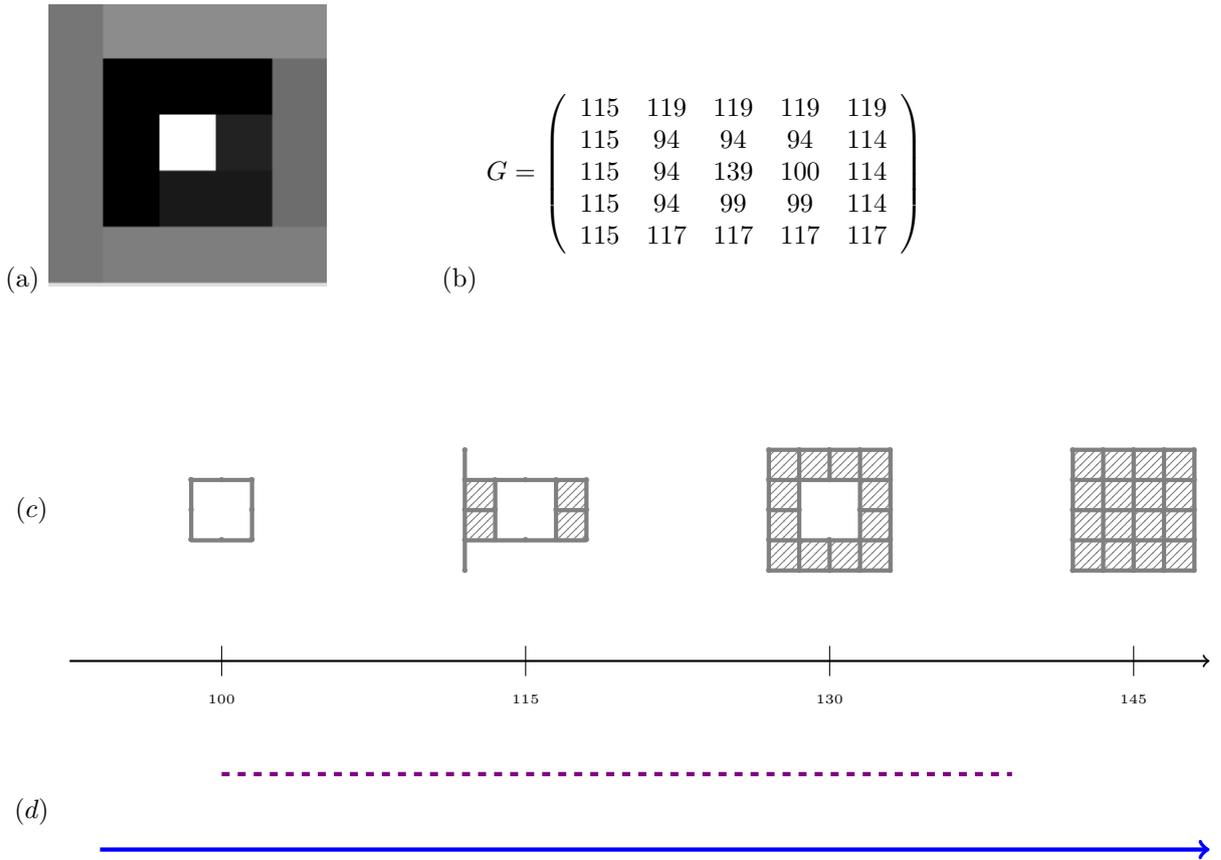
\begin{figure}
(a) \includegraphics[scale=0.7]{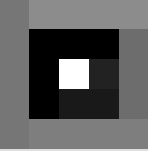}  \qquad\qquad
(b) \raisebox{4em}{$G=\left(\begin{array}{ccccc}
  115  & 119 &  119 &  119 &  119\\
   115  &  94 &   94 &   94 &  114\\
   115  &  94 &  139 &  100 &  114\\
   115   & 94   & 99   & 99  & 114\\
   115  & 117 &  117 &  117 &  117
  \end{array}\right)$}\\
  
  \vspace{1.5cm}

  \begin{tikzpicture}
\node at (-0.5,2) {$(c)$};
\node at (-0.5,-2) {$(d)$};
\draw[->, thick] (0,0)--(15,0);
\draw[-] (2,-0.2)--(2,0.2);
\draw[-] (6,-0.2)--(6,0.2);
\draw[-] (10,-0.2)--(10,0.2);  
\draw[-] (14,-0.2)--(14,0.2);  
\node at (2,-0.5) {\tiny{$100$}};
\node at (6,-0.5) {\tiny{$115$}};
\node at (10,-0.5) {\tiny{$130$}};   
\node at (14,-0.5) {\tiny{$145$}};   
\node at (2,2) {$
\begin{tikzpicture}[scale=0.4,every node/.style={scale=0.5}]
 \foreach \x in {0,1,2}
 \node at (\x,0) {\color{gray}{$\bullet$}};
  \foreach \x in {0,1,2}
 \node at (\x,2){\color{gray}{$\bullet$}};
\foreach \y in {0,1,2}
\node at (0,\y) {\color{gray}{$\bullet$}};
\foreach \y in {0,1,2}
\node at (2,\y) {\color{gray}{$\bullet$}};
\foreach \y in {0,2}
\draw[-,color=gray,ultra thick] (0,\y)--(2,\y);
\foreach \x in {0,2}
\draw[-,color=gray,ultra thick] (\x,0)--(\x,2);
\end{tikzpicture}
$};
\node at (6,2) {$
\begin{tikzpicture}[scale=0.4,every node/.style={scale=0.5}]
\foreach \y in {0,1,2,3,4}
\node at (0,\y) {\color{gray}{$\bullet$}};
\foreach \y in {1,2,3}
\node at (1,\y) {\color{gray}{$\bullet$}};
\foreach \y in {1,3}
\node at (2,\y) {\color{gray}{$\bullet$}};
\foreach \y in {1,2,3}
\node at (3,\y) {\color{gray}{$\bullet$}};
\foreach \y in {1,2,3}
\node at (4,\y) {\color{gray}{$\bullet$}};
\draw[-,color=gray,ultra thick] (0,1)--(4,1);
\draw[-,color=gray,ultra thick] (0,0)--(0,4);
\draw[-,color=gray,ultra thick] (0,3)--(4,3);
\draw[-,color=gray,ultra thick] (1,1)--(1,3);
\draw[-,color=gray,ultra thick] (3,1)--(3,3);
\draw[-,color=gray,ultra thick] (4,1)--(4,3);
\draw[-,color=gray,ultra thick] (0,2)--(1,2);
\draw[-,color=gray,ultra thick] (3,2)--(4,2);
\draw[draw=gray,pattern=north east lines, pattern color=gray, thick] (0,1) rectangle (1,3);
\draw[draw=gray,pattern=north east lines, pattern color=gray, thick] (3,1) rectangle (4,3);
\end{tikzpicture}
$};
\node at (10,2) {$
\begin{tikzpicture}[scale=0.4,every node/.style={scale=0.5}]
\foreach \y in {0,1,2,3,4}
\node at (0,\y) {\color{gray}{$\bullet$}};
\foreach \y in {0,1,2,3,4}
\node at (1,\y) {\color{gray}{$\bullet$}};
\foreach \y in {0,1,3,4}
\node at (2,\y) {\color{gray}{$\bullet$}};
\foreach \y in {0,1,2,3,4}
\node at (3,\y) {\color{gray}{$\bullet$}};
\foreach \y in {0,1,2,3,4}
\node at (4,\y) {\color{gray}{$\bullet$}};
\draw[-,color=gray,ultra thick] (0,0)--(4,0);
\draw[-,color=gray,ultra thick] (0,1)--(4,1);
\draw[-,color=gray,ultra thick] (0,3)--(4,3);
\draw[-,color=gray,ultra thick] (0,4)--(4,4);
\draw[-,color=gray,ultra thick] (0,0)--(0,4);
\draw[-,color=gray,ultra thick] (1,0)--(1,4);
\draw[-,color=gray,ultra thick] (3,0)--(3,4);
\draw[-,color=gray,ultra thick] (4,0)--(4,4);
\draw[-,color=gray,ultra thick] (0,2)--(1,2);
\draw[-,color=gray,ultra thick] (3,2)--(4,2);
\draw[-,color=gray,ultra thick] (2,0)--(2,1);
\draw[-,color=gray,ultra thick] (2,3)--(2,4);
\draw[draw=gray,pattern=north east lines, pattern color=gray, thick] (0,0) rectangle (1,4);
\draw[draw=gray,pattern=north east lines, pattern color=gray, thick] (3,0) rectangle (4,4);
\draw[draw=gray,pattern=north east lines, pattern color=gray, thick] (1,0) rectangle (3,1);
\draw[draw=gray,pattern=north east lines, pattern color=gray, thick] (1,3) rectangle (3,4);
\end{tikzpicture}
$};
\node at (14,2) {$
\begin{tikzpicture}[scale=0.4,every node/.style={scale=0.5}]
\foreach \y in {0,1,2,3,4}
\node at (0,\y) {\color{gray}{$\bullet$}};
\foreach \y in {0,1,2,3,4}
\node at (1,\y) {\color{gray}{$\bullet$}};
\foreach \y in {0,1,2,3,4}
\node at (2,\y) {\color{gray}{$\bullet$}};
\foreach \y in {0,1,2,3,4}
\node at (3,\y) {\color{gray}{$\bullet$}};
\foreach \y in {0,1,2,3,4}
\node at (4,\y) {\color{gray}{$\bullet$}};
\draw[-,color=gray,ultra thick] (0,0)--(4,0);
\draw[-,color=gray,ultra thick] (0,1)--(4,1);
\draw[-,color=gray,ultra thick] (0,3)--(4,3);
\draw[-,color=gray,ultra thick] (0,4)--(4,4);
\draw[-,color=gray,ultra thick] (0,0)--(0,4);
\draw[-,color=gray,ultra thick] (1,0)--(1,4);
\draw[-,color=gray,ultra thick] (3,0)--(3,4);
\draw[-,color=gray,ultra thick] (4,0)--(4,4);
\draw[-,color=gray,ultra thick] (0,2)--(4,2);
\draw[-,color=gray,ultra thick] (2,0)--(2,4);
\draw[draw=gray,pattern=north east lines, pattern color=gray, thick] (0,0) rectangle (4,4);

\end{tikzpicture}
$};
\draw[color=violet,ultra thick,dashed] (2,-1.5)--(12.4,-1.5) ;
\draw[->,color=blue,ultra thick] (0.4,-2.5)--(15,-2.5) ;
\end{tikzpicture}
\caption{(a) A gray-scale image, (b) the matrix of gray values, (c) the filtered cubical complex associated to the digital image, and (d) the barcode for the nested sequence of spaces in panel (c). A solid line represents the lifetime of a component, and a dashed line represents the lifetime of a hole.}\label{digital image}
\end{figure}

We now discuss a second example related to digital images. (For an illustration, see Fig.~\ref{digital image}(a).)
Digital images have a  ``cubical'' structure, given by the pixels (for $2$-dimensional digital images) or voxels (for $3$-dimensional images). Therefore, one approach to study digital images uses combinatorial structures called ``cubical complexes.'' (For a different approach to the study of digital images, see Section \ref{S:data}.) Roughly, cubical complexes are topological spaces built from a union of vertices, edges, squares, cubes, and higher-dimensional hypercubes. 
An efficient way \cite{WCV12} to build a cubical complex from a $2$-dimensional digital image consists of assigning a vertex to every pixel, then joining vertices corresponding to adjacent pixels by an edge, and filling in the resulting squares. 
One proceeds in a similar way for $3$-dimensional images. One then labels every vertex with an integer that corresponds to the gray value of the pixel, and one labels edges (respectively, squares) with the maximum of the values of the adjacent vertices (respectively, edges). 
One can then construct a nested sequence of cubical complexes  $C_0\subset C_1 \subset \dots \subset C_{256}$,
 where for each $i\in \{0,1,\dots , 256\}$, the cubical complex $C_i$ contains all vertices, edges, squares, and cubes that are labeled by a number less than or equal to $i$. (See Fig.~\ref{digital image}(c) for an example.) Such a sequence of cubical complexes is also called a ``filtered cubical complex.'' 
Similar to the previous example, one can use homology to measure several features of the spaces $C_i$ (see Fig.~\ref{digital image}(d)).

In the present article, we focus on persistent homology, but there are also other methods in TDA --- including the Mapper algorithm \cite{mapper}, Euler calculus (see \cite{G14} for an introduction with an eye towards applications), cellular sheaves \cite{G14,C13}, and many more.  We refer readers who wish to learn more about the foundations of TDA to the 
article \cite{C09}, which discusses why topology and functoriality are essential for data analysis. We point to several introductory papers, books, and two videos on PH at the end of Section \ref{S:PH}.

The first algorithm for the computation of PH was introduced for computation over $\mathbb{F}_2$ (the field with two elements) in \cite{ELZ02} and over general fields in \cite{ZC05}. Since then, several algorithms and optimization techniques have been presented, and there are now various powerful implementations of PH \cite{phat, dipha, dionysus, perseus, javaplex,gudhi,ripser}. Those wishing to try PH for computations may find it difficult  to discern which implementations and algorithms are best suited for a given task. The field of PH is evolving continually, and new software implementations and updates are released at a rapid pace.  Not all of them are well-documented, and (as is well-known  in the TDA community), the computation of PH for large data sets is computationally very expensive.

To our knowledge, there exists neither an overview of the various computational methods for PH nor a comprehensive benchmarking of the state-of-the-art implementations for the computation of persistent homology. In the present article, we close this gap: we introduce computation of PH to a  general audience of applied mathematicians and computational scientists, offer guidelines for the computation of PH, and test the existing open-source published libraries for the computation of PH.

 The rest of our paper is organized as follows. In Section \ref{S:related}, we discuss related work. We then introduce homology in Section \ref{S:homology} and introduce PH in Section \ref{S:PH}. We discuss the various steps of the pipeline for the computation of PH in Section \ref{S:computation}, and we briefly examine algorithms for ``generalized persistence''  in Section \ref{S:generalised}. In Section \ref{S:software}, we give an overview of software  libraries, discuss our benchmarking of a collection of them, and provide guidelines for which software or algorithm is better suited to which data set. (We provide specific guidelines for the computation of PH with the different libraries in the Tutorial in the Supplementary Information (SI).)  In Section \ref{S:considerations}, we discuss future directions for the computation of PH.


\section{Related work}\label{S:related}

In our work, we introduce PH to non-experts with an eye towards applications, and we benchmark state-of-the-art libraries for the computation of PH. In this section, we discuss related work for both of these points.

There are several excellent introductions to the theory of PH (see the references at the end of Section \ref{SS:filtered and homology}), but none of them emphasizes the actual computation of PH by providing specific guidelines for people who want to do computations. 
In the present paper, we navigate the theory of PH with an eye towards applications, and we provide guidelines for the computation of PH using the open-source libraries $\perseus$, $\dionysus$, $\dipha$, $\javaplex$, and $\gudhi$. We include a tutorial (see the SI) that gives specific instructions for how to use the different functionalities that are implemented in these libraries. Much of this information is scattered throughout numerous different papers, websites, and even source code of implementations, and we believe that it is beneficial to the applied mathematics community (especially people who seek an entry point into PH) to find all of this information in one place. The functionalities that we cover include plots of barcodes and persistence diagrams and the computation of PH with Vietoris--Rips complexes, alpha complexes, \v{C}ech complexes, witness complexes, and cubical complexes for image data. 
 We thus believe that our paper closes a gap in introducing PH to people interested in applications, while our tutorial complements existing tutorials (see, e.g. \cite{FKLM14,BD15,javaplex-tutorial}).

We believe that there is a need for a thorough benchmarking of the state-of-the-art libraries. In our work, we use twelve different data sets to test and compare the libraries $\javaplex$, $\gudhi$, $\dipha$, $\dionysus$, and $\perseus$, and we obtain some surprising results (see Section \ref{Conclusion}). 
There are several benchmarkings in the PH literature; we are aware of the following ones: the benchmarking in 
\cite{dSMVJ11b} compares the implementations of standard and dual algorithms in $\dionysus$; 
the one in \cite{VN12} compares the Morse-theoretic reduction algorithm with the standard algorithm;  
the one in \cite{phat} compares all of the data structures and algorithms implemented in $\phat$; 
the benchmarking in \cite{BKRJ14} compares $\phat$ and its spin-off $\dipha$; 
and the benchmarking in C.~Maria's doctoral thesis \cite{CM14} is
to our knowledge the only existing benchmarking that compares packages from different authors.
However, Maria compares only up to three different implementations at one time, and he used the package $\jplex$ (which is no longer maintained) instead of the $\javaplex$ library (its successor). Additionally, the widely used library $\perseus$ (e.g., it was used in \cite{WZDZ14,KGKM13,KGKM14,TKHKMPM14}) does not appear in Maria's benchmarking.


\section{Homology}\label{S:homology}

Assume that one is given data that lies in a metric space, such as a subset of Euclidean space with an inherited distance function. 
 In many situations, one is not interested in the precise geometry of these spaces, but instead  seeks to understand some basic characteristics, such as the number of components or the existence of holes and voids. Algebraic topology captures these basic characteristics either by counting them or by associating vector spaces or more sophisticated algebraic structures to them. Here we are interested in {\it homology}, which associates one vector space $H_i(X)$ to a space $X$ for each natural number $i \in \{0, 1,2, \dots\}$. The dimension of $H_0(X)$ counts the number of path components in $X$, the dimension of $H_1(X)$ is a count of  the number of holes, and the dimension of $H_2(X)$ is a count of the number of voids. 
An important property of these algebraic structures is that they are robust, as they do not change when the underlying space is transformed by bending, stretching, or other deformations. In technical terms, they are {\it homotopy invariant}. 
\footnote{Conversely, under favorable conditions (see \cite[Corollary 4.33]{H02}), these algebraic invariants determine the topology of a  space up to homotopy--- an equivalence relation that is much coarser (and easier to work with) than the more familiar notion of homeomorphy.}

It can be very difficult to compute the homology of arbitrary topological spaces. We thus approximate our spaces by combinatorial structures called ``simplicial complexes,'' for which homology can be easily computed algorithmically. Indeed, often one is not even given the space $X$, but instead  possesses only a discrete sample set $S$ from which to build a simplicial complex following one of the recipes described in Sections \ref{SS:complexes} and \ref{subsection: filtrations}.


\subsection{Simplicial complexes and their homology} \label{SS:simplicial complexes}

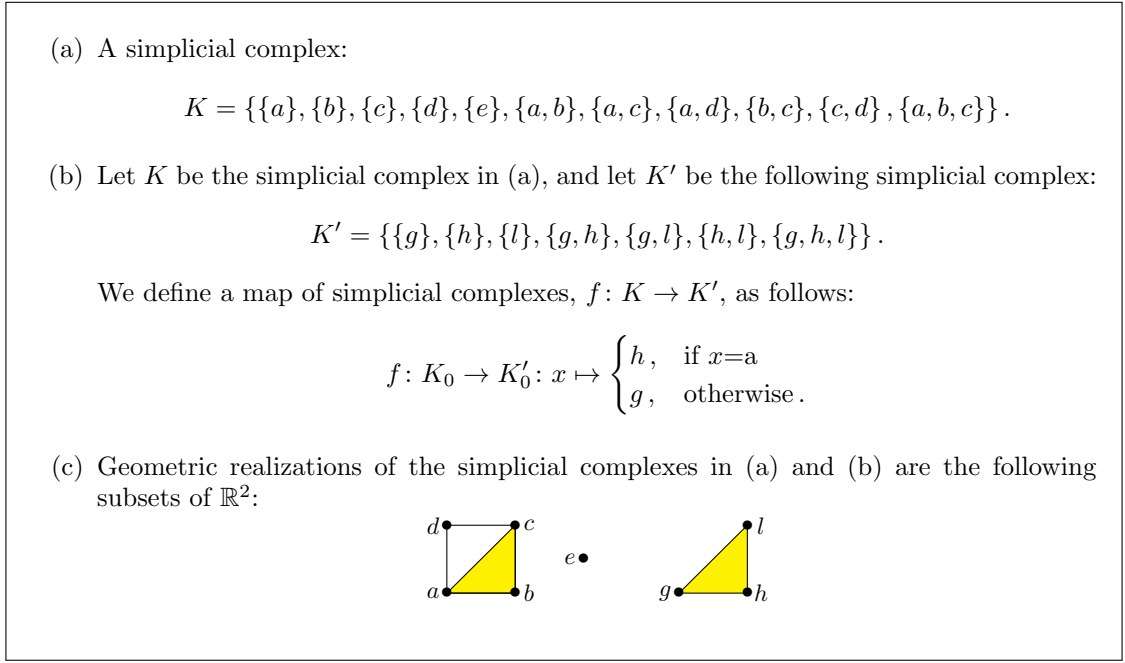
\begin{figure}[h!]
\begin{framed}
\begin{enumerate}[label=(\alph*)]
\item{\label{item: ex complex}  A simplicial complex:
\begin{equation*}
	K= \{ \{ a\}, \{ b\} , \{ c\}, \{ d\}, \{e \} , \{a,b\}, \{a, c\}, \{ a, d\} , \{b , c\}, \{c, d\}\,,
\{a,b,c\}\}\,.
\end{equation*}
}
\item{\label{item: ex map} Let $K$ be the simplicial complex in \ref{item: ex complex}, and let $K'$ be the following simplicial complex:
\begin{equation*}
	K'=\{\{ g\}, \{h\}, \{l\}, \{g,h\},\{g,l\},\{h,l\},\{g,h,l\}\}\,.
\end{equation*}
We define a map of simplicial complexes, $f\colon K\to K'$, as follows: 
\begin{equation*}
	f\colon K_0\to K'_0\colon x\mapsto \begin{cases} h\,, &\text{if $x$=a}\\
g\,, &\text{otherwise}\,.
\end{cases}
\end{equation*}
}
\item{\label{item: ex geom real} Geometric realizations of the simplicial complexes in \ref{item: ex complex} and \ref{item: ex map} are the following subsets of $\mathbb{R}^2$:
\[
\begin{tikzpicture}[scale=0.9, every node/.style={scale=0.9}]
\node at (-0.2,0) {$a$};
\node at (1.2,0) {$b$};
\node at (-0.2,1) {$d$};
\node at (1.2,1) {$c$};
\node at (1.8,0.5) {$e$};
\draw[fill=yellow](0,0)--(1,0)--(1,1)--(0,0);
\draw[] (0,0)-- (0,1)--(1,1) --(1,0)--(0,0);
\foreach \position in {(0,0), (1,0), (0,1) ,(1,1), (2,0.5)}
\node at \position{$\bullet$};
\end{tikzpicture} \qquad
\begin{tikzpicture}[scale=0.9, every node/.style={scale=0.9}]
\node at (-0.2,0) {$g$};
\node at (1.2,0) {$h$};
\node at (1.2,1) {$l$};
\draw[fill=yellow](0,0)--(1,0)--(1,1)--(0,0);
\foreach \position in {(0,0), (1,0),(1,1)}
\node at \position{$\bullet$};
\end{tikzpicture}
\] 
}
 \end{enumerate}
 \end{framed}
    \caption{(a) Example of a simplicial complex, (b) a map of simplicial complexes, and (c) a geometric realization of the simplicial complex in (a).}\label{ex: simplicial complex}
\end{figure}

We begin by giving the definitions of simplicial complexes and of the maps between them.
Roughly, a simplicial complex is a space that is built from a union of points, edges, triangles, tetrahedra, and higher-dimensional polytopes. We illustrate the main definitions given in this section with the example in Fig.~\ref{ex: simplicial complex}. As we pointed out in Section \ref{intro}, ``cubical complexes'' give another way to associate a combinatorial structure to a topological space. In TDA, cubical complexes have been used primarily to study image data sets. One can compute PH for a nested sequence of cubical complexes in a similar way as for simplicial complexes, but the theory of PH for simplicial complexes is richer, and we therefore examine only simplicial homology and complexes in our discussions. See \cite{KMM04} for a treatment of cubical complexes and their homology.

 \begin{definition}\label{def1}
A \emph{simplicial complex}\footnote{Note that  this is usually called an``abstract simplicial complex'' in the literature.} is a collection $K$ of non-empty subsets of a set $K_0$ such that $\tau\subset\sigma$ and $\sigma\in K$ guarantees that $\tau\in K$ and $\{v\}\in K$ for all $v\in K_0$. 
The elements of $K_0$ are called \emph{vertices} of $K$, and the elements of $K$ are called \emph{simplices}.  Additionally, we say that a simplex has \emph{dimension $p$} or is a \emph{$p$-simplex} if it has a cardinality of $p+1$. We use $K_p$ to denote the collection of $p$-simplices. The \emph{$k$-skeleton} of $K$ is the union of the sets $K_p$ for all $p\in \{0,1,\dots , k\}$. If $\tau$ and $\sigma$ are simplices such that $\tau\subset \sigma$, then we call $\tau$ a \emph{face} of $\sigma$. 
The \emph{dimension} of $K$ is defined as the maximum of the dimensions of its simplices.
A \emph{map of simplicial complexes}, $f : K \to L$,  is a map $f \colon K_0\to L_0$ such that $f(\sigma)\in L$ for all $\sigma\in K$.
\end{definition}

 We give an example of a simplicial complex in Fig.~\ref{ex: simplicial complex}\ref{item: ex complex} and an example of a map of simplicial complexes in Fig.~\ref{ex: simplicial complex}\ref{item: ex map}.
Definition \ref{def1} is rather abstract, but one can always interpret a finite simplicial complex $K$ geometrically as a subset of $\mathbb{R}^N$ for sufficiently large $N$; such a subset  is called a ``geometric realization,'' and it is unique up to a canonical piecewise-linear homeomorphism. 
For example, the simplicial complex in Fig.~\ref{ex: simplicial complex}\ref{item: ex complex} has a geometric realization given by the subset of $\mathbb{R}^2$ in Fig.~\ref{ex: simplicial complex}\ref{item: ex geom real}.

 We now define homology for simplicial complexes. Let $\mathbb{F}_2$ denote the field with two elements.
Given a simplicial complex $K$, let $C_p(K)$ denote the $\mathbb{F}_2$-vector space with basis given by the $p$-simplices of $K$. For any $p\in \{1,2,\dots \}$, we define the linear map
\begin{equation*}
	{d}_p\colon {C}_p(K)\to {C}_{p-1}(K)\colon \sigma\mapsto \sum_{\tau\subset \sigma, \tau \in K_{p-1}}\tau\,.
 \end{equation*}
\noindent
For  $p=0$,  we define $d_0$ to be the zero map. In words, ${d}_p$ maps each $p$-simplex to its boundary, the sum of its faces of codimension $1$. Because the boundary of a boundary is always empty, the linear maps $d_p$ have the property that composing any two consecutive maps yields the zero map: for all $p\in \{0,1,2,\dots \}$, we have ${d}_p\circ {d}_{p+1}=0$. Consequently, the image of ${d}_{p+1}$ is contained in the kernel of ${d}_p$, so we can take the quotient of $\kernel({d}_p)$ by $\image({d}_{p+1})$.  We can thus make the following definition.

 \begin{definition}
For any $p\in \{0,1,2,\dots \}$,  the \emph{$p$th  homology} of a simplicial complex $K$ is the quotient vector space 
\begin{equation*}
	H_p(K) := \kernel({d}_p)\left.\right/\image({d}_{p+1})\,.
\end{equation*}
\noindent
Its dimension 
\begin{equation*}
	\beta _p(K) :=  \dim  H_p(K) =\dim \kernel (d_p) - \dim \image (d_{p+1})
\end{equation*}
\noindent
is called the \emph{$p$th Betti number} of $K$.
Elements in the image of $d_{p+1}$ are called \emph{$p$-boundaries}, and elements in the kernel of $d_p$ are called \emph{$p$-cycles}. 
\end{definition}

 Intuitively, the $p$-cycles that are not boundaries represent $p$-dimensional holes. Therefore, the $p$th Betti number ``counts''  the number of $p$-holes. Additionally, if $K$ is a simplicial complex of dimension $n$, then for all $p>n$, we have that $H_p(K)=0$, as $K_p$ is empty and hence  $C_p (K) = 0$.
 We therefore obtain the following sequence of vector spaces and linear maps:
  \begin{equation*}
	0\overset{d_{n+1}}{\longrightarrow}C_n(K)\overset{d_n}{\longrightarrow}\dots \dots \overset{d_2}{\longrightarrow}C_1(K)\overset{d_1}{\longrightarrow}C_0(K)\overset{d_0}{\longrightarrow} 0\;.
 \end{equation*}

\noindent
We give an example of such a sequence in Fig.~\ref{figure: simpl homology}\ref{item: ex hom}, for which we also report the Betti numbers.

\begin{figure}
[h!]
\begin{framed}
\begin{enumerate}[label=(\alph*)]
\item  \label{item: ex hom} We  compute the simplicial homology for the simplicial complex in Fig.~\ref{ex: simplicial complex}\ref{item: ex complex}. We have the following sequence of vector spaces and linear maps: 
\begin{equation*}
	0\longrightarrow \mathbb{F}_2\overset{d_2}{\longrightarrow} \mathbb{F}_2^5\overset{d_1}{\longrightarrow}  \mathbb{F}_2^5 \overset{d_0}{\longrightarrow} 0\,. 
\end{equation*}
Let $abc$ denote the basis vector that corresponds to the simplex $\{a,b,c\}$. Similarly, we use $ab$, $ac$, $ad$, $bc$, and $cd$ to  denote the basis vectors that correspond to the $1$-simplices; and we use $a$, $b$, $c$, $d$, and $e$ to denote the basis vectors that correspond to the $0$-simplices. We order the bases of the vector spaces using lexicographic order. We then have
\begin{equation*}
	d_2= \begin{pmatrix}
  1&1&0&1&0
 	\end{pmatrix}^\transpose
 \end{equation*}
  and 
  \begin{equation*}
 	 d_1= \begin{pmatrix}
  1&1&1&0&0\\
  1 &0&0&1&0\\
  0&1&0&1&1\\
  0&0&1&0&1\\
 	 0&0&0&0&0
   \end{pmatrix}\,.
 \end{equation*} 
 One can then compute that $\beta _0 (K) = 2$, $\beta _1(K) =1$, and all higher Betti numbers are $0$.
 \item \label{ex:hom map} The map of simplicial complexes, $f\colon K\to K'$, in Fig.~\ref{ex: simplicial complex}\ref{item: ex map} induces the map $\widetilde{f_0}\colon C_0(K)\to C_0(K')$ on chain complexes that sends the basis element $a$ to $h$ (note that we use the same notation as in part \ref{item: ex hom} of this box). Furthermore, it sends the basis elements that correspond to the other $0$-simplices of $K$ to $g$, and it induces the map $\widetilde{f_1}\colon C_1(K)\to C_1(K')$ that sends the basis elements $ab$, $ac$, and $ad$ to $cd$ and sends the basis elements $cd$ and $bc$ to $0$. This gives a map $ f_{0}$ on degree-$0$ homology 
 that sends both generators of $H_0(K)$ to the generator of $H_0(K')$. It follows that $H_i(K')=0$ for all $i\geq 1$, so $f_i$ is the zero map for $i\geq 1$. 
\end{enumerate}
\end{framed}
\caption{(a) Computation of simplicial homology for the simplicial complex in Fig.~\ref{ex: simplicial complex}\ref{item: ex complex} and (b) induced map in degree-$0$ homology 
for the map of simplicial complexes in Fig.~\ref{ex: simplicial complex}\ref{item: ex map}.}\label{figure: simpl homology}
\end{figure}

One of the most important properties of simplicial homology is ``functoriality.''  Any map $f\colon K\to K'$ of simplicial complexes induces the following $\mathbb{F}_2$-linear map:
\begin{equation*}
	\widetilde{f}_p\colon  C_p(K)\to C_p(K')\colon \sum_{\sigma\in K_p} c_\sigma \sigma \mapsto \sum_{\sigma\in K_p \text{ such that } f(\sigma)\in K'_p} c_\sigma f(\sigma) \,, \quad \text{for any} \quad p\in \{0,1,2,\dots\} \,,
\end{equation*}
where $c_\sigma\in \mathbb{F}_2$. Additionally, $\widetilde{f}_{p}\circ d_{p+1}=d'_{p+1}\circ \widetilde{f}_{p+1}$, and the map $\widetilde{f}_p$ therefore induces the following linear map between homology vector spaces:
\begin{align*}
	f_{ p}\colon H_p(K)&\to H_p(K') \,, \notag \\
	[c] &\mapsto [\widetilde{f}_p(c)] \;.
\end{align*} 
(We give an example of such a map in Fig.~\ref{figure: simpl homology}\ref{ex:hom map}.) Consequently, to any map $f\colon K\to K'$ of simplicial complexes, we can assign a map $f_p\colon H_p(K)\to H_p(K')$ for any $p\in \{0,1,2,\dots\}$. This assignment has the important property that given a pair of composable maps of simplicial complexes, $f\colon K \longrightarrow K'$ and $g\colon K' \longrightarrow K''$, the map $(g \circ f)_p\colon H_p(K)\to H_p(K'')$ is equal to the composition of the maps induced by $f$ and $g$. That is, $(g \circ f)_p=g_p\circ f_p$.  The fact that a map of simplicial complexes induces a map on homology that is compatible with composition is called \emph{functoriality}, and it is crucial for the definition of persistent homology (see Section \ref{SS:filtered and homology}).

When working with simplicial complexes, one can modify a simplicial complex 
by removing or adding a pair of simplices $(\sigma, \tau)$,  where $\tau $ is a codimension-1 face of $\sigma$ and $\sigma$ is the only simplex that has $\tau$ as a face.  The resulting simplicial complex has the same homology as the one with which we started. In Fig.~\ref{ex: simplicial complex} (a), we can remove the pair $(\{a, b, c\}, \{ b,c\})$ and then the pair $(\{ a,b\}, \{b \})$ without changing  the Betti numbers. Such a move is called an {\it elementary simplicial collapse} \cite{cohen}. In Section \ref{SS:reduction}, we will see an application of this for the computation of PH. 

In this section, we have defined simplicial homology over the field $\mathbb{F}_2$ --- i.e., ``with coefficients in  $\mathbb{F}_2$.'' One can be more general and instead define simplicial homology with coefficients in any field (or even in the integers). However, when $1 \neq -1$, one needs to take more care when defining the boundary maps $d_p$ to ensure that $d_p \circ d_{p+1}$ remains the zero map. Consequently, the definition is more involved. 
For the purposes of the present paper, it suffices to consider homology with coefficients in the field $\mathbb{F}_2$. Indeed, we will see in Section \ref{S:PH} that to obtain topological summaries in the form of barcodes, we need to compute homology with coefficients in a field. Furthermore, as we summarize in Table \ref{table: software} (in Section \ref{S:software}), 
most of the implementations for the computation of PH work with $\mathbb{F}_2$. 

We conclude this section with a warning: changing the coefficient field can affect the Betti numbers. For example, if one computes the homology of the Klein bottle (see Section \ref{S:data sets}) with coefficients in the field $\mathbb{F}_p$ with $p$ elements, where $p$ is a prime, then $\beta_0(K)=1$ for all primes $p$. However, $\beta_1(K)=2$ and $\beta_2(K)=1$ if $p=2$, but $\beta_1(K)=1$ and $\beta_2(K)=0$ for all other primes $p$.  The fact that $\beta_2(K)=0$ for $p\ne 2$ arises from the nonorientability of the Klein bottle. 
 The treatment of 
different coefficient fields is beyond the scope of our article, but interested readers can peruse \cite{H02} for an introduction to homology and \cite{KMM04} for an overview of computational homology.


\subsection{Building simplicial complexes}\label{SS:complexes}

As we discussed in Section \ref{SS:simplicial complexes}, computing the homology of finite simplicial complexes boils down to linear algebra. The same is not true for the homology of an arbitrary space $X$, and one therefore tries to find simplicial complexes whose homology approximates the homology of the space in an appropriate sense.

An important tool is the \v{C}ech (\v{C}) complex. Let $\collection$ be a cover of $X$ --- i.e., a collection of subsets of $X$ such that the union of the subsets is $X$. The $k$-simplices of the \v{C}ech complex are the non-empty intersections of $k+1$ sets in the cover $\collection$. More precisely, we define the nerve of a collection of sets as follows.
\begin{definition}
Let $\collection=\{U_i\}_{i\in I}$ be a non-empty collection of sets. The \emph{nerve} of $\collection$ is the  simplicial complex with set of vertices given by $I$ and $k$-simplices given by $\{i_0,\dots , i_k\}$ if and only if $U_{i_0}\cap \dots \cap U_{i_k}\ne \emptyset$.
\end{definition}

If the cover of the sets is sufficiently ``nice,'' then the Nerve Theorem implies that the nerve of the cover and the space $X$ have the same homology \cite{B95, EH10}. For example, suppose that we have a finite set of points $S$ in a metric space $X$. We then can define, for every $\epsilon>0$, the space $S_\epsilon$ as the union $\cup_{x\in S}B(x,\epsilon)$, where $B(x,\epsilon)$ denotes the closed ball with radius $\epsilon$ centered at $x$. It follows that $\{B(x,\epsilon)\mid x\in S\}$ is a cover of $S_\epsilon$, and the nerve of this cover is the {\it \v{C}ech complex on $S$ at scale $\epsilon$}. We denote this complex by $\cech_\epsilon(S)$. 
If the space $X$ is Euclidean space, then the Nerve Theorem guarantees that the simplicial complex $\cech_\epsilon(S)$ recovers the homology of $S_\epsilon$. 

From a computational point of view, the \v{C}ech complex is expensive because one has to check for large numbers of intersections. Additionally, in the worst case, the \v{C}ech complex can have dimension $|\collection|-1$, and it therefore can have many simplices in dimensions higher than the dimension of the underlying space. Ideally, it is desirable to construct simplicial complexes that approximate the homology of a space but are easy to compute and have ``few'' simplices, especially in high dimensions. This is a subject of ongoing research:
In Subsection \ref{subsection: filtrations}, we give an overview of state-of-the-art methods to associate complexes to point-cloud data in a way that addresses one or both of these desiderata. See \cite{ES52,EH10} for more details on the \v{C}ech complex, and see \cite{B95, EH10} for a precise statement of the Nerve Theorem.


\section{Persistent homology}\label{S:PH}

Assume that we are given experimental data in the form of a finite metric space $S$; there are points or vectors that represent measurements along with some distance function (e.g., given by a correlation or a measure of dissimilarity) on the set of points or vectors. Whether or not the set $S$ is a sample from some underlying topological space, it is useful to think of it in those terms. Our goal is to recover the properties of such an underlying space in a way that is robust to small perturbations in the data $S$. In a broad sense, this is the subject of {\it topological inference}. (See\cite{Ou15} for an overview.) If $S$ is a subset of Euclidean space, one can consider a ``thickening'' $S_\epsilon$ of $S$ given by the union of balls of a certain fixed radius $\epsilon$ around its points and then compute the \v{C}ech complex. 
One can thus try to compute qualitative features of the data set $S$ by constructing the \v{C}ech complex for a chosen value $\epsilon$ and then computing its simplicial homology. The problem with this approach is that there is a priori no clear choice for the value of the parameter $\epsilon$. The key insight of PH is the following: To extract qualitative information from data, one considers several (or even all) possible values of the parameter $\epsilon$. As the value of $\epsilon$ increases, simplices are added to the complexes. Persistent homology then captures how the homology of the complexes changes as the parameter value increases, and it detects which features ``persist'' across changes in the parameter value. We give an example of persistent homology in Fig.~\ref{figure: PH}.


\subsection{Filtered complexes and homology}\label{SS:filtered and homology}

Let $K$ be a finite simplicial complex, and let $K_1 \subset K_2 \subset \dots \subset K_l= K$ be a finite sequence of
nested subcomplexes of $K$. The simplicial complex $K$ with such a sequence of subcomplexes is called a {\it filtered simplicial complex}. See Fig.~\ref{figure: PH}(a) for an example of filtered simplicial complex.
We can apply homology to each of the subcomplexes. For all homology degrees $p$, the inclusion maps $K_i\to K_j$ induce $\mathbb{F}_2$-linear maps $f_{i,j}\colon H_p(K_i)\to H_p(K_j)$ for all $i,j \in \{1,\dots,l\}$ with $i \leq j$. By functoriality (see Section \ref{SS:simplicial complexes}), it follows that 
\begin{equation}\label{eq:pers mod}
	f_{k,j}\circ f_{i,k}=f_{i,j} \qquad \text{for all $i\leq k\leq j$\,.}
\end{equation} 
We therefore give the following definition\footnote{A pair $\left( \{M_i\}_{i\in I}, \{\phi_{i,j}\colon M_i\to M_j\}_{i\leq j}\right)$, where $(I,\leq )$ is a totally ordered set, such that for each $i$, we have that $M_i$ is a vector space and the maps $\phi_{i,j}$ are linear maps satisfying the functoriality property (\ref{eq:pers mod}), is usually called a \emph{persistence module}. With this terminology, the homology of a filtered simplicial complex is an example of persistence module.}.

\begin{definition}
Let $K_1 \subset K_2 \subset \dots \subset K_l= K$ be a filtered simplicial complex. The \emph{$p$th persistent homology of $K$} is the pair 
\begin{equation*}
	\left(\left\{H_p(K_i)\right\}_{1\leq i\leq l},\{f_{i,j}\}_{1\leq i\leq j\leq l}\right)\,,
\end{equation*} 
where for all $i,j \in \{1,\dots , l\}$ with $i\leq j$, the linear maps $f_{i,j}\colon H_p(K_i)\to H_p(K_j)$ are the maps induced by the inclusion maps $K_i\to K_j$.
\end{definition}

\begin{center}
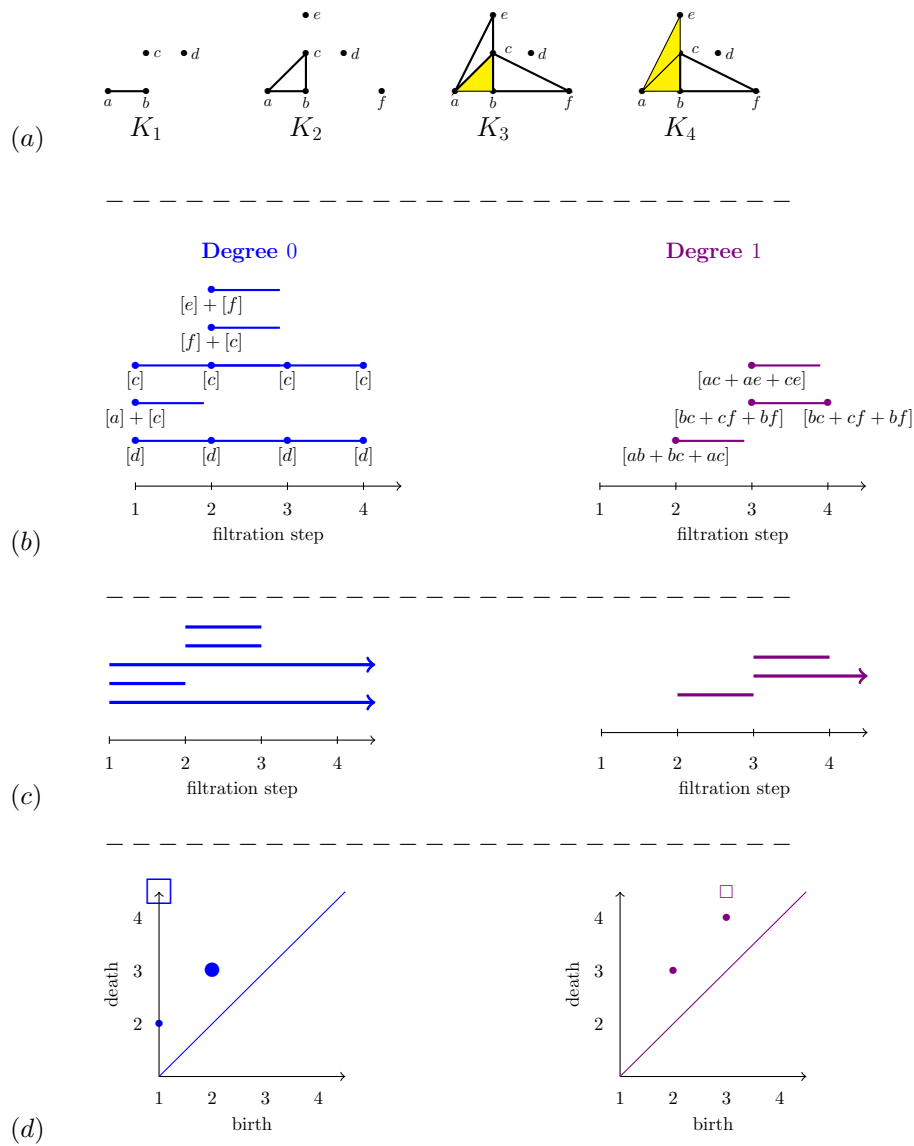
\begin{figure}[htbp!]
\begin{alignat}{-1}
&\notag (a)\qquad
&&\notag
\begin{tikzpicture}[scale=0.5, every node/.style={scale=0.6}]
\foreach \position in {(1,1),(2,1),(2,2),(3,2)}
\node at \position{$\bullet$};
\node at (1,0.7) {$a$};
\node at (2,0.7) {$b$};
\node at (2.3,2) {$c$};
\node at (3.3,2) {$d$};
\draw[thick] (1,1)--(2,1);
\node at (2,0) {\LARGE{$K_1$}};
\end{tikzpicture} 
\qquad
\begin{tikzpicture}[scale=0.5, every node/.style={scale=0.6}]
\foreach \position in {(1,1),(2,1),(2,2),(3,2),(4,1),(2,3)}
\node at \position{$\bullet$};
\node at (1,0.7) {$a$};
\node at (2,0.7) {$b$};
\node at (2.3,2) {$c$};
\node at (3.3,2) {$d$};
\node at (4,0.7) {$f$};
\node at (2.3,3) {$e$};
\node at (2,0) {\LARGE{$K_2$}};
\draw[thick] (1,1)--(2,1)--(2,2)--(1,1);
\end{tikzpicture}
\qquad
\begin{tikzpicture}[scale=0.5, every node/.style={scale=0.6}]
\foreach \position in {(1,1),(2,1),(2,2),(3,2),(4,1),(2,3)}
\node at \position{$\bullet$};
\node at (1,0.7) {$a$};
\node at (2,0.7) {$b$};
\node at (2.4,2.1) {$c$};
\node at (3.3,2) {$d$};
\node at (4,0.7) {$f$};
\node at (2.3,3) {$e$};
\node at (2,0) {\LARGE{$K_3$}};
\draw[fill=yellow](1,1)--(2,1)--(2,2)--(1,1);
\draw[thick] (2,3)--(2,2)--(1,1)--(2,3);
\draw[thick](2,1)--(4,1)--(2,2)--(2,1);
\end{tikzpicture} 
\qquad
\begin{tikzpicture}[scale=0.5, every node/.style={scale=0.6}]
\foreach \position in {(1,1),(2,1),(2,2),(3,2),(4,1),(2,3)}
\node at \position{$\bullet$};
\node at (1,0.7) {$a$};
\node at (2,0.7) {$b$};
\node at (2.4,2.1) {$c$};
\node at (3.3,2) {$d$};
\node at (4,0.7) {$f$};
\node at (2.3,3) {$e$};
\node at (2,0) {\LARGE{$K_4$}};
\draw[fill=yellow](1,1)--(2,1)--(2,2)--(1,1);
\draw[fill=yellow](2,3)--(2,2)--(1,1)--(2,3);
\draw[thick](2,1)--(4,1)--(2,2)--(2,1);
\end{tikzpicture}\\
&\notag &&\notag -------------------------- \\
&\notag (b) \qquad
&&\notag
\begin{tikzpicture}[scale=0.5, every node/.style={scale=0.7}]
\node at (3,6){\color{blue}{\textbf{\large{Degree $0$}}}};
\node at (0,-0.8){1};
\node at (2,-0.8){2};
\node at (4,-0.8){3};
\node at (6,-0.8){4};
\draw[->] (0,-0.2)--(7,-0.2);
\draw[-] (0,-0.3)--(0,-0.1);
\draw[-] (2,-0.3)--(2,-0.1);
\draw[-] (4,-0.3)--(4,-0.1);
\draw[-] (6,-0.3)--(6,-0.1);
\node at (3.5,-1.5) {filtration step};
\node at (0,1) {\color{blue}{$\bullet$}};
\node at (0,0.6) {$[d]$};
\node at (0,2) {\color{blue}{$\bullet$}};
\node at (0,1.6) {$[a]+[c]$};
\node at (0,3) {\color{blue}{$\bullet$}};
\node at (0,2.6) {$[c]$};
\node at (2,1) {\color{blue}{$\bullet$}};
\node at (2,0.6) {$[d]$};
\node at (2,3) {\color{blue}{$\bullet$}};
\node at (2,2.6) {$[c]$};
\node at (2,4) {\color{blue}{$\bullet$}};
\node at (2,3.6) {$[f]+[c]$};
\node at (2,5) {\color{blue}{$\bullet$}};
\node at (2,4.6) {$[e]+[f]$};
\node at (4,3) {\color{blue}{$\bullet$}};
\node at (4,2.6) {$[c]$};
\node at (4,1) {\color{blue}{$\bullet$}};
\node at (4,0.6) {$[d]$};
\node at (6,3) {\color{blue}{$\bullet$}};
\node at (6,2.6) {$[c]$};
\node at (6,1) {\color{blue}{$\bullet$}};
\node at (6,0.6) {$[d]$};
\path[-, thick, color=blue] (0,1) edge (2,1)
(2,1) edge (4,1)
(4,1) edge (6,1)
(0,3) edge  (2,3)
(0,2) edge  (1.8,2)
(2,3) edge(4,3)
(4,3) edge (6,3)
(2,3) edge (3.8,3)
(2,4) edge (3.8,4)
(2,5) edge (3.8,5);
\end{tikzpicture} \qquad\qquad\qquad\quad
\begin{tikzpicture}[scale=0.5, every node/.style={scale=0.7}]
\node at (3,6){\color{violet}{\textbf{\large{Degree $1$}}}};
\node at (0,-0.8){1};
\node at (2,-0.8){2};
\node at (4,-0.8){3};
\node at (6,-0.8){4};
\draw[->] (0,-0.2)--(7,-0.2);
\draw[-] (0,-0.3)--(0,-0.1);
\draw[-] (2,-0.3)--(2,-0.1);
\draw[-] (4,-0.3)--(4,-0.1);
\draw[-] (6,-0.3)--(6,-0.1);
\node at (3.5,-1.5) {filtration step};
\node at (2,1) {\color{violet}{$\bullet$}};
\node at (2,0.6) {$[ab+bc+ac]$};
\node at (4,2) {\color{violet}{$\bullet$}};
\node at (3.4,1.6) {$[bc+cf+bf]$};
\node at (4,3) {\color{violet}{$\bullet$}};
\node at (4,2.6) {$[ac+ae+ce]$};
\node at (6,2) {\color{violet}{$\bullet$}};
\node at (6.8,1.6) {$[bc+cf+bf]$};
\path[-,color=violet, thick] 
(2,1) edge (3.8,1)
(4,2) edge (6,2)
(4,3) edge (5.8,3);
\end{tikzpicture}
 \\
&\notag &&\notag -------------------------- \\
&\notag (c) \qquad
&&\notag
\begin{tikzpicture}[scale=0.5, every node/.style={scale=0.7}]
\node at (0,-0.6){1};
\node at (2,-0.6){2};
\node at (4,-0.6){3};
\node at (6,-0.6){4};
\draw[->] (0,0)--(7,0);
\draw[-] (0,-0.1)--(0,0.1);
\draw[-] (2,-0.1)--(2,0.1);
\draw[-] (4,-0.1)--(4,0.1);
\draw[-] (6,-0.1)--(6,0.1);
\node at (3.5,-1.3) {filtration step};
\draw[->,color=blue,very thick](0,1)--(7,1);
\draw[color=blue,very thick](0,1.5)--(2,1.5);
\draw[->,color=blue,very thick] (0,2)--(7,2);
\draw[color=blue,very thick](2,2.5)--(4,2.5);
\draw[color=blue,very thick](2,3)--(4,3);
 \end{tikzpicture}  
 \qquad\qquad\qquad\qquad
  \begin{tikzpicture}[scale=0.5, every node/.style={scale=0.7}]
\node at (0,-0.8){1};
\node at (2,-0.8){2};
\node at (4,-0.8){3};
\node at (6,-0.8){4};
\draw[->] (0,0)--(7,0);
\draw[-] (0,-0.1)--(0,0.1);
\draw[-] (2,-0.1)--(2,0.1);
\draw[-] (4,-0.1)--(4,0.1);
\draw[-] (6,-0.1)--(6,0.1);
\node at (3.5,-1.5) {filtration step};
 \draw[color=violet, very thick] (2,1)--(4,1) ;
   \draw[->,color=violet, very thick] (4,1.5)--(7,1.5) ;
  \draw[color=violet,very thick] (4,2)--(6,2) ;
 \end{tikzpicture} \\
&\notag &&\notag -------------------------- \\
&\notag (d) 
&&\notag\begin{tikzpicture}[scale=0.35, every node/.style={scale=0.7}]
\node at (0,-0.8){1};
\node at (2,-0.8){2};
\node at (4,-0.8){3};
\node at (6,-0.8){4};
\node at (-0.8,2){2};
\node at (-0.8,4){3};
\node at (-0.8,6){4};
\node at (0,2) {\color{blue}{$\bullet$}};
\node at (0,7) {\color{blue}{\scalebox{2}{$\Box$}}};
\node at (2,4) {\color{blue}{\scalebox{2}{$\bullet$}}};
\draw[->] (0,0)--(7,0);
 \draw[->] (0,0)--(0,7) ;
 \draw[-,color=blue] (0,0)--(7,7);
\node at (3.5,-1.8) {birth};
\node at (-1.8,3.5) {\rotatebox{90}{death}};
\node at (3,4){};
 \end{tikzpicture} \qquad\qquad\qquad\qquad
\begin{tikzpicture}[scale=0.35, every node/.style={scale=0.7}]
\node at (0,-0.8){1};
\node at (2,-0.8){2};
\node at (4,-0.8){3};
\node at (6,-0.8){4};
\node at (-0.8,2){2};
\node at (-0.8,4){3};
\node at (-0.8,6){4};
\node at (2,4){\color{violet}{$\bullet$}};
\node at (4,6) {\color{violet}{$\bullet$}};
\node at (4,7) {\color{violet}{$\Box$}};
 \draw[->] (0,0)--(7,0) ;
  \draw[->] (0,0)--(0,7) ;
  \draw[-,color=violet] (0,0)--(7,7);
 \node at (3.5,-1.8) {birth};
 \node at (-1.8,3.5) {\rotatebox{90}{death}};
  \node at (3,4){};
 \end{tikzpicture}
 \end{alignat}
  \caption{Example of persistent homology for a finite filtered simplicial complex. (a) We start with a finite filtered simplicial complex. 
  (b) At each filtration step $i$, we draw as many vertices as the dimension of (left column) $H_0(K_i)$ and (right column) $H_1(K_i)$. We label the vertices by basis elements, the existence of which is guaranteed by the Fundamental Theorem of Persistent Homology, and  we draw an edge between two vertices to represent the maps $f_{i,j}$, as explained in the main text. We thus obtain a well-defined collection of disjoint half-open intervals called a ``barcode.'' We interpret each interval in degree $p$ as representing the lifetime of a $p$-homology class across the filtration.   
  (c) We rewrite the diagram in (b) in the conventional way.  We represent classes that are born but do not die at the final filtration step using arrows that start at the birth of that feature and point to the right. 
  (d) An alternative graphical way to represent barcodes (which gives exactly the same information) is to use {\it persistence diagrams}, in which an interval $[i,j)$ is represented by the point $(i,j)$ in the extended plane $\overline{\mathbb{R}}^2$, where $\overline{\mathbb{R}}=\mathbb{R}\cup \{\infty\}$. Therefore, a persistence diagram is a finite multiset of points in $\overline{\mathbb{R}}^2$. We use squares to signify the classes that do not die at the final step of a filtration, and the size of dots or squares is directly proportional to the number of points being represented. For technical reasons, which we discuss briefly in Section \ref{SS:interpretation}, one also adds points on the diagonal to the persistence diagrams. (Each of the points on the diagonal has infinite multiplicity.)} \label{figure: PH}
 \end{figure} 
  \end{center}

The $p$th persistent homology of a filtered simplicial complex gives more refined information than just the homology of the single subcomplexes. We can visualize the information given by the vector spaces $H_p(K_i)$ together with the linear maps $f_{i,j}$ by drawing the following diagram: at filtration step $i$, we draw as many bullets as the dimension of the vector space $H_p(K_i)$. We then connect the bullets as follows: we draw an interval between bullet  $u$ at filtration step $i$ and bullet  $v$ at filtration step $i+1$ if the generator of $H_p(K_i)$ that corresponds to $u$ is sent to the generator of $H_p(K_{i+1})$ that corresponds to $v$. If the generator corresponding to a bullet $u$ at filtration step $i$ is sent to $0$ by $f_{i,i+1}$, we draw an interval starting at  $u$ and  ending at $i+1$.
(See Fig.~\ref{figure: PH}(b) for an example.) Such a diagram clearly depends on a choice of basis for the vector spaces $H_p(K_i)$, and a poor choice can lead to complicated and unreadable clutter. Fortunately, by the Fundamental Theorem of Persistent Homology \cite{ZC05}, there is a choice of basis vectors of $H_p(K_i)$ for each $i\in \{1,\dots , l\}$ such that 
one can construct the diagram as a well-defined and unique collection of disjoint half-open intervals, collectively called a {\it barcode}\footnote{Although the collection of intervals is unique, note that one has to choose a vertical order when drawing the intervals in the diagram, and there is therefore an ambiguity in the representation of the intervals as a barcode. However, there is no ambiguity when representing the intervals as points in a persistence diagram (see Fig.~\ref{figure: PH}(d)).}.  We give an example of a barcode in Fig.~\ref{figure: PH}(c). Note that the Fundamental Theorem of PH, and hence the existence of a barcode, relies on the fact that we are using homology with field coefficients. (See \cite{ZC05} for more details.) 

There is a useful interpretation of barcodes in terms of births and deaths of generators.
Considering the maps $f_{i,j}$ written in the basis given by the Fundamental Theorem of Persistent Homology, we say that $x\in H_p(K_i)$ (with $x \neq 0$) is \emph{born} in $H_p(K_i)$ if it is not in the image of $f_{i-1,i}$ (i.e., $f^{-1}_{i-1,i}(x)=\emptyset$). For $x \in H_p(K_i)$ (with $x \neq 0$), we say that $x$ \emph{dies} in $H_p(K_j)$ if $j>i$ is the smallest index for which $f_{i,j}(x)=0$. The lifetime of $x$ is represented by the half-open interval $[i,j)$. If $f_{i,j}(x)\ne 0$ for all $j>i$ in $I$, we say that $x$ \emph{lives forever}, and its lifetime is represented by the interval $[i,\infty)$.

\begin{rmk}
Note that some references (e.g., \cite{EH10}) introduce persistent homology by defining  the birth and death of generators without using the existence of a choice of compatible bases, as given by the Fundamental Theorem of Persistent Homology. The definition of birth coincides with the definition that we have given, but the definition of death is different. One says that $x \in H_p(K_i)$ (with $x \neq 0$) \emph{dies} in $H_p(K_j)$ if $j>i$ is the smallest index for which either $f_{i,j}(x)=0$ or there exists $y\in H_p(K_{i'})$ with $i'<i$ such that $f_{i',j}(y)=f_{i,j}(x)$. In words, this means that $x$ and $y$ merge at filtration step $j$, and the class that was born earlier is the one that survives. In the literature, this is called the \emph{elder rule}. We do not adopt this definition, because the elder rule is not well-defined when two classes are born at the same time, as there is no way to choose which class will survive. For example, in Fig.~\ref{figure: PH}, there are two classes in $H_0$ that are born at the same stage in $K_1$. These two classes merge in $K_2$, but neither dies. The class that dies is $[a]+[c]$. 
\end{rmk}

There are numerous excellent introductions to PH, such as the books \cite{EH10,G14,Z09,Ou15} and the papers \cite{weinberger,G08,EH08,EM12,C09}.
 For a brief and friendly introduction to PH and some of its applications, see the video \url{https://www.youtube.com/watch?v=h0bnG1Wavag}. For a brief introduction to some of the ideas in TDA, see the video \url{https://www.youtube.com/watch?v=XfWibrh6stw}.

\section{Computation of PH for data}\label{S:computation}

\begin{figure}[h!]
\begin{center}
\begin{tikzpicture}[scale=1, every node/.style={scale=1}][transform shape]
\node at (0,0) {}; 
\path \etape{1}{ Data};
\path (p1.east)+(2.0,0.0)   \etape{2}{ Filtered complex};
\path [arrow] (p1.east) -- node [above] {$(1)$} (p2);
\path (p2.east)+(2.0,0.0) \etape{3}{ Barcodes};
\path (p3.east)+(2.2,0.0) \etapee{4}{Interpretation};
\path [arrow] (p2.east) -- node [above] {$(2)$} (p3);
\path [arrow] (p3.east) -- node [above] {$(3)$} (p4);
\end{tikzpicture}
\end{center}
\caption{PH pipeline.}\label{figure:first pipeline}
\end{figure}
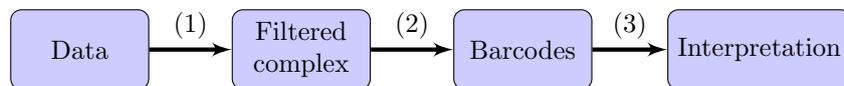

We summarize the pipeline for the computation of PH from data in Fig.~\ref{figure:first pipeline}. In the following subsections, we describe each step of this pipeline and state-of-the-art algorithms for the computation of PH. The two features that make PH appealing for applications are that it is computable via linear algebra 
and that it is stable with respect to perturbations in the measurement of data. 
In Section \ref{SS:stability}, we give a brief overview of stability results.

\subsection{Data}\label{S:data}

As we mentioned in Section \ref{intro}, types of data sets that one can study with PH include finite metric spaces, digital images, and networks. We now give a brief overview of how one can study these types of data sets using PH.  


\subsubsection{Networks}

One can construe an undirected network as a $1$-dimensional simplicial complex. If the network is weighted, then filtering by increasing or decreasing weight yields a filtered $1$-dimensional simplicial complex. To obtain more refined information about the network, it is desirable to construct higher-dimensional simplices. There are various methods to do this. The simplest method, called a {\it weight rank clique filtration} (WRCF), consists of building a clique complex on each subnetwork. (See Section \ref{SSS:VR} for the definition of ``clique complex.'') See \cite{PSDV13} for an application of this method.  Another method to study networks with PH consists of mapping the nodes of the network to points of a finite metric space. There are several ways to compute distances between nodes of a network; the method that we use in our benchmarking in Section \ref{S:software} consists of computing a shortest path between nodes. 
For such a distance to be well-defined, note that one needs the network to be connected (although conventionally one takes the distance between nodes in different components to be infinity).
There are many methods to associate an unfiltered simplicial complex to both undirected and directed networks. See the book \cite{JJ07} for an overview of such methods, and see the paper \cite{HMR09} for an overview of PH for networks.


\subsubsection{Digital images}

As we mentioned in Section \ref{intro}, digital images have a natural cubical structure: $2$-dimensional digital images are made of pixels, and $3$-dimensional images are made of voxels. Therefore, to study digital images, cubical complexes are more appropriate than simplicial complexes. Roughly, cubical complexes are spaces built from a union of vertices, edges, squares, cubes, and so on. 
One can compute PH for cubical complexes in a similar way as for simplicial complexes, and we will therefore not discuss this further in this paper. See \cite{KMM04} for a treatment of computational homology with cubical complexes rather than simplicial complexes and for a discussion of the relationship between simplicial and cubical homology. See \cite{WCV12} for an efficient algorithm and data structure for the computation of PH for cubical data, and \cite{BEK10} for an algorithm that computes PH for cubical data in an approximate way. For an application of PH and cubical complexes to movies, see \cite{VN12}.

Other approaches for studying digital images are also useful. In general, given a digital image that consists of $N$ pixels or voxels, one can consider this image as a point in a $c\times N$-dimensional space, with each coordinate storing a vector of length $c$ representing the color of a pixel or voxel. Defining an appropriate distance function on such a space allows one to consider a collection of images (each of which has $N$ pixels or voxels) as a finite metric space. A  version of this approach was used in \cite{CIdSZ08}, in which the local structure of natural images was studied by selecting $3\times 3$ patches of pixels of the images. 


\subsubsection{Finite metric spaces}\label{SSS:finite metric spaces}

As we mentioned in the previous two subsections, both undirected networks and image data can be construed as finite metric spaces. Therefore, methods to study finite metric spaces with PH apply to the study of networks and image data sets.

In some applications, points of a metric space have associated ``weights.'' For instance, in the study of molecules, one can represent a molecule as a union of balls in Euclidean space \cite{ZY12,xia16}.
For such data sets, one would therefore also consider a minimum filtration value (see Section \ref{subsection: filtrations} for the description of such filtration values) at which the point enters the filtration. In Table \ref{table: software}(g), we indicate which software libraries implement this feature.


\subsection{Filtered Simplicial Complexes}\label{subsection: filtrations}

In Section \ref{SS:complexes}, we introduced the \v{C}ech complex, a classical simplicial complex from algebraic topology.
However, there are many other simplicial complexes that are better suited for studying data from applications. We  discuss them in this section.

To be a useful tool for the study of data, a simplicial complex has to satisfy some theoretical properties dictated by topological inference; roughly, if we build the simplicial complex on a set of points sampled from a space, then the homology of the simplicial complex  has to approximate the homology of the space. 
 For the \v{C}ech complex, these properties are guaranteed by the Nerve Theorem. Some of the complexes that we discuss in this subsection are motivated by a ``sparsification paradigm'': they approximate the PH of known simplicial complexes but have fewer simplices than them. Others, like the Vietoris--Rips complex, are appealing because they can be computed efficiently. 
In this subsection, we also review reduction techniques, which are heuristics that reduce the size of complexes without changing the PH. In Table \ref{table:complexes}, we summarize the simplicial complexes that we discuss in this subsection.

For the rest of this subsection $(X,\dist)$ denotes a metric space, and $S$ is a subset of $X$, which becomes a metric space with the induced metric. In applications, $S$ is the collection of measurements together with a notion of distance, and we assume that $S$ lies in the (unknown) metric space $X$. Our goal is then to compute persistent homology for a sequence of nested spaces $S_{\epsilon_1},S_{\epsilon_2},\dots , S_{\epsilon_l}$, where each space gives a ``thickening'' of $S$ in $X$.

\subsubsection{Vietoris--Rips complex}\label{SSS:VR}

We have seen that one of the disadvantages of the \v{C}ech complex is that one has to check for a large number of intersections. To circumvent this issue, one can instead consider the Vietoris--Rips (VR) complex, which approximates the \v{C}ech complex. For a non-negative real number $\epsilon$, the {\it Vietoris--Rips complex $\VR_\epsilon(S)$ at scale $\epsilon$} is defined as
\begin{equation*}
	\VR_\epsilon(S)=\left\{\sigma\subseteq S \mid \dist(x,y)\leq 2\epsilon \text{ for all }x,y \in \sigma\right\}\,.
\end{equation*}

\noindent

The sense in which the VR complex approximates the \v{C}ech complex is that, when $S$ is a subset of Euclidean space, we have 
$\cech_\epsilon (S) \subseteq \VR_\epsilon(S) \subseteq \cech_{\sqrt {2}\epsilon} (S)$.
Deciding whether a subset $\sigma\subseteq S$ is in $\VR_\epsilon(S)$ is equivalent to deciding if the maximal pairwise distance between any two vertices in $\sigma$ is at most $2\epsilon$. Therefore, one can construct the VR complex in two steps. One first computes the {\it $\epsilon$-neighborhood graph of $S$}. This is the graph whose vertices are all points in $S$ and whose edges are
\begin{equation*}
	\{(i,j)\in S\times S\mid i\ne j \text{ and } \dist(i,j)\leq 2\epsilon\}\,.
\end{equation*}
Second, one obtains the VR complex by computing the {\it clique complex} of the $\epsilon$-neighborhood graph. The clique complex of a graph is a simplicial complex that is defined as follows: The subset $\{x_0,\dots , x_k\}$ is a $k$-simplex if and only if every pair of vertices in $\{x_0,\dots , x_k\}$ is connected by an edge. Such a collection of vertices is called a \emph{clique}. This construction makes it very easy to compute the VR complex, because to construct the clique complex one has only to check for pairwise distances --- for this reason clique complexes are also called ``lazy'' in the literature. Unfortunately, the VR complex has the same worst-case complexity as the \v{C}ech complex. In the worst case, it can have up to $2^{|S|}-1$ simplices and dimension $|S|-1$. 

In applications, one therefore usually only computes the VR complex up to some dimension $k \ll |S|-1$. In our benchmarking, we often choose $k=2$ and $k=3$.

The paper \cite{Z10} overviews different algorithms to perform both of the steps for the construction of the VR complex, and it introduces fast algorithms to construct the clique complex. For more details on the VR complex, see \cite{V27,EH10}. For a proof of the approximation of the \v{C}ech complex by the VR complex, see \cite{EH10}; see \cite{KS13} for a generalization of this result.


\subsubsection{The Delaunay complex}

To avoid the computational problems of the \v{C}ech and VR complexes, we need a way to limit the number of simplices in high dimensions. The Delaunay complex gives a geometric tool to accomplish this task, and most of the new simplicial complexes that have been introduced for the study of data are based on variations of the Delaunay complex. The Delaunay complex and its dual, the Voronoi diagram, are central objects of study in computational geometry because they have many useful properties.

For the Delaunay complex, one usually considers $X=\mathbb{R}^d$, so we also make this assumption. We subdivide the space $\mathbb{R}^d$ into regions of points that are closest to any of the points in $S$. More precisely, for any $s\in S$, we define 
\begin{equation*}
	V_s=\{x\in \mathbb{R}^d \mid \dist(x,s)\leq \dist(x,s') \text{ for all } s'\in S \}\,.
\end{equation*}
The collection of sets $V_s$ is a cover for $\mathbb{R}^d$ that is called the {\it Voronoi decomposition of $\mathbb{R}^d$ with respect to $S$}, and the nerve of this cover is called the {\it Delaunay complex} of $S$ and is denoted by $\del(S;\mathbb{R}^d)$. In general, the Delaunay complex does not have a geometric realization in $\mathbb{R}^d$. However, if the points $S$ are ``in general position''
\footnote{A set $S$ of points in $\mathbb{R}^d$ is \emph{in general position} if no $d+2$ points of $S$ lie on a $d$-dimensional sphere, and for any $d'<d$, no $d'+2$ points of $S$ lie on a $d'$-dimensional subspace that is isometric to $\mathbb{R}^{d'}$. In particular, a set of points $S$ in $\mathbb{R}^2$ is in general position if no four points lie on a $2$-dimensional sphere and no three points lie on a line.} 
then the Delaunay complex has a geometric realization in $\mathbb{R}^d$ that gives a triangulation of the convex hull of $S$. In this case, the Delaunay complex is also called the {\it Delaunay triangulation}.

The complexity of the Delaunay complex depends on the dimension $d$ of the space. For $d\leq 2$, the best algorithms have complexity $\bigo(N \log N)$, where $N$ is the cardinality of $S$. For $d\geq 3$, they have complexity $\bigo(N^{\lceil d/2 \rceil})$. The construction of the Delaunay complex is therefore costly in high dimensions, although there are efficient algorithms for the computation of the Delaunay complex for $d = 2$ and $d = 3$. Developing efficient algorithms for the construction of the Delaunay complex in higher dimensions is a subject of ongoing research. See \cite{BDH09} for a discussion of progress in this direction, and see \cite{GOR97} for more details on the Delaunay complex and the Voronoi diagram.


\subsubsection{Alpha complex}\label{SS:alpha}

We continue to assume that $S$ is a finite set of points in $\mathbb{R}^d$. Using the Voronoi decomposition, one can define a simplicial complex that is similar to the \v{C}ech complex, but which has the desired property that (if the points $S$ are in general position) its dimension is at most that of the space. Let $\epsilon>0$, and let $S_\epsilon$ denote the union $\bigcup_{s\in S}B(s,\epsilon)$. For every $s\in S$, consider the intersection $V_s\cap B(s,\epsilon)$. The collection of these sets forms a cover of $S_\epsilon$, and the nerve complex of this cover is called the {\it alpha ($\alpha$) complex of $S$ at scale $\epsilon$} and is denoted by $A_\epsilon(S)$. The Nerve Theorem applies, and it therefore follows that $A_\epsilon(S)$ has the same homology as $S_\epsilon$. 

Furthermore, $A_\infty(S)$ is the Delaunay complex; and for $\epsilon<\infty$, the alpha complex is a subcomplex of the Delaunay complex. The alpha complex was introduced for points in the plane in \cite{EKS83}, in 3-dimensional Euclidean space in \cite{EM94}, and for Euclidean spaces of arbitrary dimension in \cite{E95}. For points in the plane, there is a well-known speed-up for the alpha complex that uses a duality between $0$-dimensional and $1$-dimensional persistence for alpha complexes \cite{EH08}. (See \cite{Kur15} for the algorithm, and see \cite{hopes} for an implementation.)


\subsubsection{Witness complexes}\label{S:witness}

Witness complexes are very useful for analyzing large data sets, because they make it possible to construct a simplicial complex on a significantly smaller subset $L\subseteq S$ of points that are called ``landmark'' points. Meanwhile, because one uses information about all points in $S$ to construct the simplicial complex, the points in $S$ are called ``witnesses.'' 
Witness complexes can be construed as a ``weak version'' of Delaunay complexes. (See the characterization of the Delaunay complex in \cite{dS08}.) 

\begin{definition} Let $(S,\dist)$ be a metric space, and let $L\subseteq S$ be a finite subset. Suppose that $\sigma$ is a non-empty subset of $L$. We then say that $s\in S$ is a \emph{weak witness for $\sigma$ with respect to $L$} if and only if $\dist(s,a)\leq \dist(s,b)$ for all $a\in \sigma$ and for all $b\in L\setminus \sigma$.
The \emph{weak Delaunay complex $\del^w(L;S)$ of $S$ with respect to $L$} has vertex set given by the points in $L$, and a subset $\sigma$ of $L$ is in $\del^w(L;S)$ if and only if it has a weak witness in $S$. 
\end{definition}

To obtain nested complexes, one can extend the definition of witnesses to $\epsilon$-witnesses.
\begin{definition}
A point $s\in S$ is a \emph{weak $\epsilon$-witness} for $\sigma$ with respect to $L$ if and only if $\dist(s,a)\leq \dist(s,b)+\epsilon$ for all $a\in \sigma$ and for all $b\in L\setminus \sigma$. \end{definition}

Now we can define the {\it weak Delaunay complex $\del^w(L;S,\epsilon)$ at scale $\epsilon$} to be the simplicial complex with vertex set $L$, and such that a subset $\sigma\subseteq L$ is in $\del^w(L;S,\epsilon)$ if and only if it has a weak $\epsilon$-witness in $S$. 
By considering different values for the parameter $\epsilon$, we thereby obtain nested simplicial complexes.
The weak Delaunay complex is also called the ``weak witness complex'' or just the ``witness complex'' in the literature. 

There is a modification of the witness complex called the {\it lazy witness complex $\del^w_{\mathrm{lazy}}(L;X,\epsilon)$}. It is a clique complex, and it can therefore be computed more efficiently than the witness complex. 
The lazy witness complex has the same $1$-skeleton as $\del^w(L;X,\epsilon)$, and one adds a simplex $\sigma$ to $\del^w_{\mathrm{lazy}}(L;X,\epsilon)$ whenever its edges are in $\del^w_{\mathrm{lazy}}(L;X,\epsilon)$. Another type of modification of the witness complex yields {\it parametrized witness complexes}. 
Let $\nu =1,2,\dots$ and for all $s\in S$ define $m_\nu(s)$ to be the distance to the $\nu$th closest landmark point.
Furthermore, define $m_0(s)=0$ for all $s\in S$. Let $\W_\nu(L;S, \epsilon)$ be the simplicial complex whose vertex set is $L$ and such that a $1$-simplex $\sigma=\{x_0,x_1\}$ is in $\W_\nu(L;X, \epsilon)$ if and only if there exists  $s$ in $S$ for which
\begin{equation*}
	\max\{\dist(x_0,s),\dist(x_1,s)\}\leq m_\nu(s)+\epsilon \,.
\end{equation*}
A simplex $\sigma$ is in $\W_\nu(L;X, \epsilon)$ if and only if all of its edges belong to $\W_\nu(L;X, \epsilon)$.
 For $\nu=2$, note that $\W_2(L;X, \epsilon)=\del^w_{\mathrm{lazy}}(L;X,\epsilon)$. For $\nu=0$, we have that $\W_0(L;X, \epsilon)$ approximates the VR complex $\VR(L;\epsilon)$. That is,
 \begin{equation*}
	\W_0(L;X,\epsilon)\subseteq \VR(L;2\epsilon)\subseteq \W_0(L;X, 2\epsilon)\,.
 \end{equation*}
 \noindent
Note that parametrized witness complexes are often called ``lazy witness complexes'' in the literature, because they are clique complexes. 

The weak Delaunay complex was introduced in \cite{dS08}, and parametrized witness complexes were introduced in \cite{dSC04}. Witness complexes can be rather useful for applications. Because their complexity depends on the number of landmark points, one can reduce complexity by computing simplicial complexes using a smaller number of vertices. 
However, there are theoretical guarantees for the witness complex only when $S$ is the metric space associated to a low-dimensional Euclidean submanifold. It has been shown that witness complexes can be used to recover the topology of curves and surfaces in Euclidean space \cite{GO08,AEM07}, but they can fail to recover topology for submanifolds of Euclidean space of 
three or more dimensions \cite{BGO09}. Consequently, there have been studies of simplicial complexes that are similar to the witness complexes but with better theoretical guarantees (see  Section \ref{SSS:add compl}).


 \subsubsection{Additional complexes}\label{SSS:add compl}
 
Many more complexes have been introduced for the fast computation of PH for large data sets. These include the graph-induced complex \cite{DTW13}, which is a simplicial complex constructed on a subsample $Q$, and has better theoretical guarantees than the witness complex (see \cite{gic} for the companion software);
an approximation of the VR complex that has a worst-case size that is linear in the number of data points \cite{S13};  an approximation of the \v{C}ech complex\cite{KS13} whose worst-case size also scales linearly in the data; and an approximation of the VR complex via simplicial collapses \cite{DSW16}. We do not discuss such complexes in detail, because thus far (at the time of writing) none of them have been implemented in publicly-available libraries for the computation of PH. (See Table \ref{table: software} in Section \ref{S:software} for information about which complexes have been implemented.) 

\begin{table}
\begin{center}
\footnotesize
\caption{We summarize several types of complexes that are used for PH. We indicate the theoretical guarantees and the worst-case sizes of the complexes as functions of the cardinality $N$ of the vertex set. For the witness complexes (see Section \ref{S:witness}), $L$ denotes the set of landmark points, while $Q$ denotes the subsample set for the graph-induced complex (see Section \ref{SSS:add compl}).
}\label{table:complexes}
\begin{tabular}{|c|c|c|}
\toprule
Complex $K$& Size of $K$ & Theoretical guarantee\\\hline
\v{C}ech &$2^{\bigo(N)}$ & Nerve theorem  \\\hline
 Vietoris--Rips (VR)  & $2^{\bigo(N)}$& Approximates \v{C}ech complex\\\hline
Alpha & \begin{tabular}{@{}c@{}}  $N^{\bigo(\lceil d/2 \rceil)}$\\ ($N$ points in $\mathbb{R}^d$)\end{tabular} & Nerve theorem \\\hline
Witness &$2^{\bigo(|L|)}$&For curves and surfaces in Euclidean space\\\hline
Graph-induced complex & $2^{\bigo(|Q|)}$& Approximates $\VR$ complex\\\hline
Sparsified \v{C}ech &$\bigo(N)$& Approximates  \v{C}ech complex\\\hline
Sparsified VR &$\bigo(N)$& Approximates  $\VR$ complex\\
\bottomrule
\end{tabular}
\end{center}
\end{table}
\normalsize


\subsubsection{Reduction techniques}\label{SS:reduction}

Thus far, we have discussed techniques to build simplicial complexes with possibly ``few'' simplices. One can also take an alternative approach to speed up the computation of PH. For example, one can use a heuristic (i.e., a method without theoretical guarantees on the speed-up) to reduce the size of a filtered complex while leaving the PH unchanged.

For simplicial complexes, one such method is based on discrete Morse theory \cite{Fo98}, which was adapted to filtrations of simplicial complexes in \cite{MN13}. The basic idea of the algorithm developed in \cite{MN13} is that one can compute a partial matching of the simplices in a filtered simplicial complex so that (i) pairs occur only between simplices that enter the filtration at the same step, (ii) unpaired simplices determine the homology, and (iii) one can remove paired simplices from the filtered complex without altering the total PH. Such deletions are examples of the elementary simplicial collapses that we mentioned in Section \ref{SS:simplicial complexes}. Unfortunately, the problem of finding an optimal partial matching was shown to be NP complete \cite{JP06}, and one thus relies on heuristics to find partial matchings to reduce the size of the complex.

One particular family of elementary collapses, called \emph{strong collapses}, was introduced in \cite{BM12b}. Strong collapses preserve cycles of shortest length in the representative class of a generator of a hole \cite{WMSK13}; this feature makes strong collapses useful for finding holes in networks \cite{WMSK13}. 
A distributed version of the algorithm introduced in \cite{WMSK13} was introduced in \cite{WCKMS13} and adapted for the computation of PH in \cite{WCK14}.

A method for the reduction of the size of a complex for clique complexes, such as the VR complex, was introduced in \cite{Z10b} and is called the {\it tidy-set method}. Using maximal cliques, this method extracts a minimal representation of the graph that determines the clique complex. Although the tidy-set method cannot be extended to filtered complexes, it can be used for the computation of zigzag PH (see Section \ref{S:generalised}) \cite{Z12}. The tidy-set method is a heuristic, because it does not give a guarantee to minimize  the size of the output complex.


\subsection{From a Filtered Simplicial Complex to Barcodes}\label{S: intervals}

To compute the PH of a filtered simplicial complex $K$ and obtain a barcode like the one illustrated in Fig.~\ref{figure: PH}(c), we need to associate to it a matrix --- the so-called {\it boundary matrix} --- that stores information about the faces of every simplex. To do this, we place a total ordering on the simplices of the complex that is compatible with the filtration in the following sense:
\begin{itemize}
\item a face of a simplex precedes the simplex;
\item a simplex in the $i$th complex $K_i$ precedes simplices in $K_j$ for $j>i$, which are not in $K_i$.
\end{itemize}
Let $n$ denote the total number of simplices in the complex, and let $\sigma_1,\dots , \sigma_n$ denote the simplices with respect to this ordering. We construct a square matrix $\delta$ of dimension $n \times n$ by storing a $1$ in $\delta(i,j)$ if the simplex $\sigma_i$ is a face of simplex $\sigma_j$ of codimension $1$; otherwise, we store a $0$ in $\delta(i,j)$.

Once one has constructed the boundary matrix, one has to reduce it using Gaussian elimination.\footnote{As we mentioned in Section \ref{S:PH}, for the reduction of the boundary matrix and thus the computation of PH, it is crucial that one uses simplicial homology with coefficients in a field; see \cite{ZC05} for details.}
In the following subsections, we discuss several algorithms for reducing the boundary matrix.


\subsubsection{Standard algorithm}\label{standard}

The so-called {\tt standard algorithm} for the computation of PH was introduced for the field $\mathbb{F}_2$ in \cite{ELZ02} and for general fields in \cite{ZC05}. For every $j \in \{1,\dots, n\}$, we define $\low(j)$ to be the largest index value $i$ such that $\delta(i,j)$ is different from $0$.\footnote{This map is called ``low'' in the literature, because one can think of it as indicating  the index of the ``lowest'' row --- the one that is nearest to the bottom of the page on which one writes the boundary matrix --- that contains a $1$ in column $j$.} If column $j$ only contains $0$ entries, then the value of $\low(j)$ is undefined. We say that the boundary matrix is {\it reduced} if the map $\low$ is injective on its domain of definition. In Alg.~\ref{figure: algorithm}, we illustrate the standard algorithm for reducing the boundary matrix. Because this algorithm operates on columns of the matrix from left to right, it is also sometimes called the ``column algorithm.'' In the worst case, the complexity of the standard algorithm is cubic in the number of simplices.

\begin{algorithm}
\begin{center}
\begin{minipage}{200pt}
\footnotesize
\begin{algorithmic}
\FOR {$j=1$ to $n$} 
\vspace{-0.3cm}
\STATE 
\WHILE {there exists $i < j$ with $\low(i) = \low(j)$}
\STATE add column $i$ to column $j$
 \ENDWHILE
\ENDFOR
\end{algorithmic}
\end{minipage}
\end{center}
\caption{The standard algorithm for the reduction of the boundary matrix to barcodes.}\label{figure: algorithm}
\end{algorithm}



\subsubsection{Reading off the intervals}

\begin{figure}[htbp!]
\begin{framed}
\begin{enumerate}[label=(\alph*)]
\item \label{item: filtered complex}
A filtered simplicial complex:
\begin{multicols}{4}
\raisebox{-1.4cm}{\begin{tikzpicture}
\node at (0,0) {$\bullet$};
\node at (0.5,-0.3) {$K_1$};
\node at (2,0.5) {$\subset$};
\end{tikzpicture}}

\columnbreak
\begin{tikzpicture}
\node at (0,0) {$\bullet$};
\node at (1,0) {$\bullet$};
\node at (0.5,1) {$\bullet$};
\node at (0.5,-0.3) {$K_2$};
\node at (2,0.5) {$\subset$};
\draw (0,0)--(0.5,1);
\end{tikzpicture}

\columnbreak
\begin{tikzpicture}
\node at (0,0) {$\bullet$};
\node at (1,0) {$\bullet$};
\node at (0.5,1) {$\bullet$};
\node at (0.5,-0.3) {$K_3$};
\node at (2,0.5) {$\subset$};
\draw (0,0)--(0.5,1);
\draw (0,0)--(1,0)--(0.5,1);
\end{tikzpicture}

\columnbreak
\begin{tikzpicture}
\node at (0,0) {$\bullet$};
\node at (1,0) {$\bullet$};
\node at (0.5,1) {$\bullet$};
\node at (0.5,-0.3) {$K_4$};
\draw (0,0)--(0.5,1);
\draw (0,0)--(1,0)--(0.5,1);
\draw[fill=yellow] (0,0)--(1,0)--(0.5,1)--(0,0);
\end{tikzpicture}
\end{multicols}

\item \label{item: order on simplices} We put a total order on the simplices that is compatible with the filtration:
\begin{multicols}{4}
\raisebox{-1.5cm}{\begin{tikzpicture}
\node at (0,0) {$\bullet$};
\node at (-0.2,-0.2) {$\sigma_1$};
\end{tikzpicture}}

\columnbreak
\begin{tikzpicture}
\node at (0,0) {$\bullet$};
\node at (1,0) {$\bullet$};
\node at (-0.2,-0.2) {$\sigma_1$};
\node at (0.5,1.3) {$\sigma_2$};
\node at (0.5,1) {$\bullet$};
\node at (1.2,-0.2) {$\sigma_3$};
\path[-]
 (0,0) edge node[left] {$\sigma_4$}(0.5,1);
\end{tikzpicture}

\columnbreak
\begin{tikzpicture}
\node at (0,0) {$\bullet$};
\node at (1,0) {$\bullet$};
\node at (-0.2,-0.2) {$\sigma_1$};
\node at (0.5,1.3) {$\sigma_2$};
\node at (0.5,1) {$\bullet$};
\node at (1.2,-0.2) {$\sigma_3$};
\path[-]
 (0,0)edge node[left]{$\sigma_4$}(0.5,1)
 (0,0) edge node[below]{$\sigma_5$}(1,0)
 (1,0) edge node[right] {$\sigma_6$}(0.5,1);
\end{tikzpicture}

\columnbreak
\begin{tikzpicture}
\node at (0,0) {$\bullet$};
\node at (1,0) {$\bullet$};
\node at (-0.2,-0.2) {$\sigma_1$};
\node at (0.5,1.3) {$\sigma_2$};
\node at (0.5,1) {$\bullet$};
\node at (1.2,-0.2) {$\sigma_3$};
\draw[fill=yellow] (0,0)--(1,0)--(0.5,1)--(0,0);
\path[-]
 (0,0)edge node[left]{$\sigma_4$}(0.5,1)
 (0,0) edge node[below]{$\sigma_5$}(1,0)
 (1,0) edge node[right] {$\sigma_6$}(0.5,1);
\node at (0.5,0.5) {$\sigma_7$};
\end{tikzpicture}
\end{multicols}
where $\sigma_i$ denotes the $i$th simplex in this order.\\

\item \label{item: boundary matrix} (Left) The boundary matrix $B$ for the filtered simplicial complex in \ref{item: filtered complex} with respect to order on simplices in \ref{item: order on simplices}, and (right) its reduction $ \overline{B}$ given by applying Alg.~\ref{figure: algorithm} (one first adds column $5$ to column $6$, and then column $4$ to column $6$):
\begin{multicols}{2}
$
 B= \begin{pmatrix}
  0&0&0&1&1&0&0\\
  0&0&0&1&0&1&0\\
  0&0&0&0&1&1&0\\
  0&0&0&0&0&0&1\\
  0&0&0&0&0&0&1\\
   0&0&0&0&0&0&1\\
   0&0&0&0&0&0&0 
   \end{pmatrix}\;
   $
   \columnbreak
   $
   \overline{B}= \begin{pmatrix}
  0&0&0&1&1&0&0\\
  0&0&0&1&0&0&0\\
  0&0&0&0&1&0&0\\
  0&0&0&0&0&0&1\\
  0&0&0&0&0&0&1\\
   0&0&0&0&0&0&1\\
   0&0&0&0&0&0&0 
   \end{pmatrix}\,.
   $
   \end{multicols}
   \item \label{item: intervals} We read off the following intervals from the matrix $\overline{B}$ in \ref{item: boundary matrix}:
\begin{itemize}
\item $\sigma_1$ is positive, unpaired; this gives the interval $[1,\infty)$ in $H_0$.
\item $\sigma_2$ is positive, paired with $\sigma_4$; this gives no interval, because $\sigma_2$ and $\sigma_4$ enter at the same time in the filtration.
\item $\sigma_3$ is positive, paired with $\sigma_5$: this gives the interval $[2,3)$ in $H_0$.
\item $\sigma_6$ is positive, paired with $\sigma_7$: this gives the interval $[3,4)$ in $H_1$.
\end{itemize}
\end{enumerate}
\end{framed}
\caption{Example of PH computation with the {\tt standard algorithm} (see Alg.~\ref{figure: algorithm}).}\label{ex: PH computation}
\end{figure}
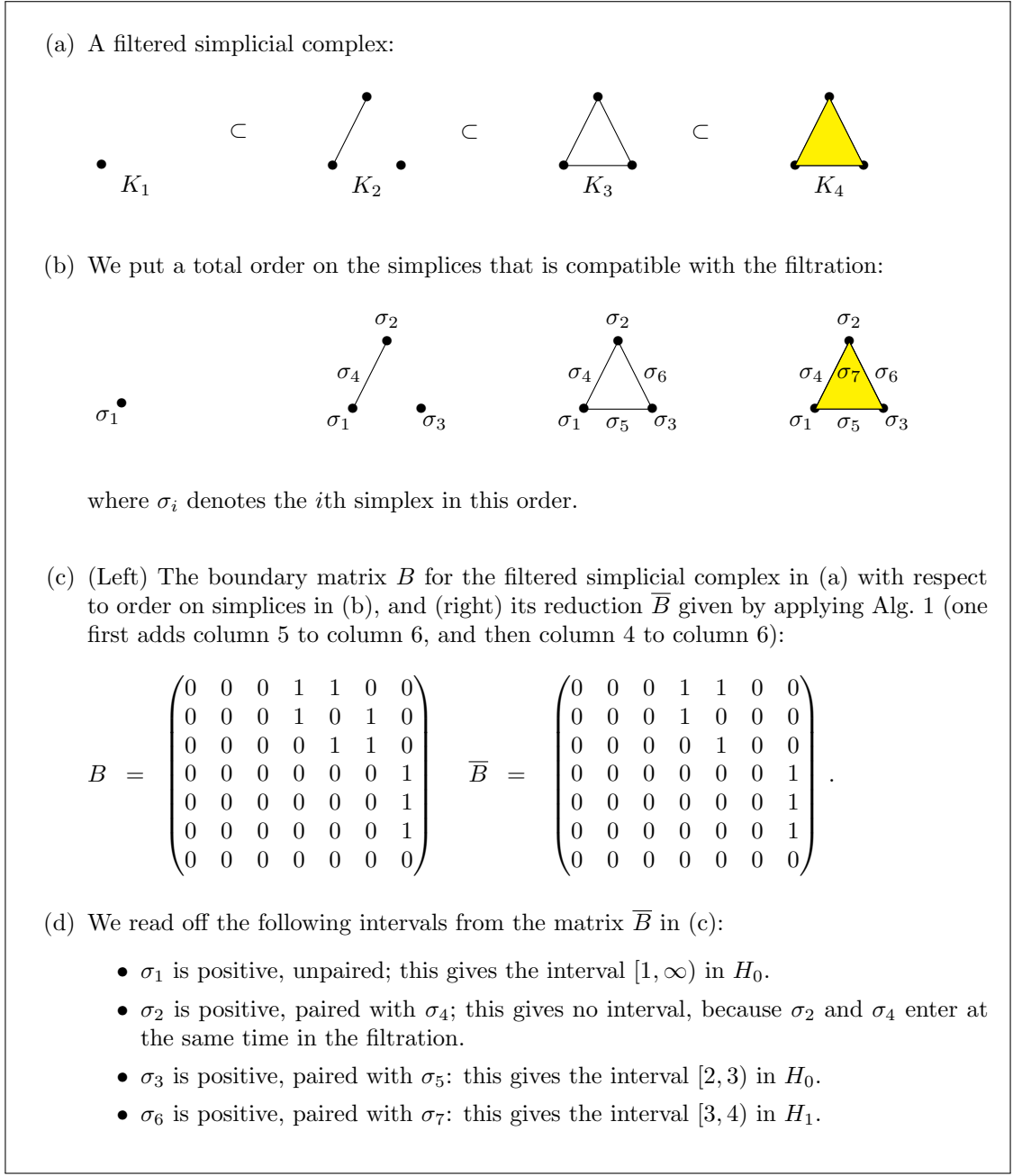

Once the boundary matrix is reduced, one can read off the intervals of the barcode by pairing the simplices in the following way:
 \begin{itemize} 
\item If $\low(j)=i$, then the simplex $\sigma_j$ is paired with $\sigma_i$, and the entrance of $\sigma_i$ in the filtration causes the birth of a feature that dies with the entrance of $\sigma_j$.
\item If $\low(j)$ is undefined, then the entrance of the simplex $\sigma_j$ in the filtration causes the birth of a feature. It there exists $k$ such that $\low(k)=j$, then $\sigma_j$ is paired with the simplex $\sigma_k$, whose entrance in the filtration causes the death of the feature. If no such $k$ exists, then $\sigma_j$ is unpaired.
 \end{itemize} 
A pair $(\sigma_i,\sigma_j)$ gives the half-open interval $[\dg(\sigma_i),\dg(\sigma_j))$ in the barcode, where for a simplex $\sigma\in K$ we define $\dg(\sigma)$ to be the smallest
number $l$ such that $\sigma\in K_l$. An unpaired simplex $\sigma_k$ gives the infinite interval $[\dg(\sigma_k),\infty)$. We give an example of PH computation in Fig.~\ref{ex: PH computation}.

\subsubsection{Other algorithms}
After the introduction of the standard algorithm, several new algorithms were developed. Each of these algorithms gives the same output for the computation of PH, so we only give a brief overview and references to these algorithms, as one does not need to know them to compute PH with one of the publicly-available software packages. In Section \ref{Conclusion}, we indicate which implementation of these libraries is best suited to which data set.

As we mentioned in Section \ref{standard}, in the worst case, the standard algorithm has cubic complexity in the number of simplices. This bound is sharp, as Morozov gave an example of a complex with cubic complexity in \cite{M05}. Note that in cases such as when matrices are sparse, complexity is less than cubic.
Milosavljevi\'c, Morozov, and Skraba \cite{MMS11} introduced an algorithm for the reduction of the boundary matrix in $\bigo(n^\omega)$, where $\omega$ is the matrix-multiplication coefficient (i.e., $\bigo(n^\omega)$ is the complexity of the multiplication of two square matrices of size $n$). 
At present, the best bound for $\omega$ is $2.376$ \cite{CW90}. Many other algorithms have been proposed for the reduction of the boundary matrix. These algorithms give a heuristic speed-up for many data sets and complexes (see the benchmarkings in the forthcoming references), but they still have cubic complexity in the number of simplices. 
Sequential algorithms include the  {\tt twist algorithm} \cite{CK11} and the {\tt dual algorithm} \cite{dSMVJ11,dSMVJ11b}. (Note that the dual algorithm is known to give a speed-up when one computes PH with the VR complex, but not necessarily for other types of complexes (see also the results of our benchmarking for the $\vertebra$ data set in the SI).) Parallel algorithms in a shared setting include the {\tt spectral-sequence algorithm} \cite[Section VII.4]{EH10} and the {\tt chunk algorithm} \cite{BKR14}, and parallel algorithms in a distributed setting include the {\tt distributed algorithm} \cite{BKRJ14}.
The {\tt multifield algorithm} is a sequential algorithm that allows the simultaneous computation of PH over several fields \cite{BM14}.


\subsection{Statistical interpretation of topological summaries}\label{SS:interpretation}

Once one has obtained barcodes, one needs to interpret the results of computations. 
In applications, one often wants to compare the output of a computation for a certain data set with the output for a null model. Alternatively, one may be studying data sets from the output of a generative model (e.g., many realizations from a model of random networks), and it is then necessary to average results over multiple realizations. In the first instance, one needs both a way to compare the two different outputs and a way to evaluate the significance of the result for the original data set. In the second case, one needs a way to calculate appropriate averages (e.g., summary statistics) of the result of the computations. 

From a statistical perspective, one can interpret a barcode as an unknown quantity that one tries to estimate by computing PH.
If one wants to use PH in applications, one thus needs a reliable way to apply statistical methods to the output of the computation of PH. To our knowledge, statistical methods for PH were addressed for the first time in the paper \cite{BK07}. 
Roughly speaking, there are three current approaches to the problem of statistical analysis of barcodes. In the first approach, researchers 
study topological properties of random simplicial complexes (see, e.g., \cite{ABW14,bianconi16}) and the review papers \cite{ABBSW10,kahle14}. One can view random simplicial complexes as null models to compare with empirical data when studying PH.  In the second approach, one studies properties of a metric space whose points are persistence diagrams. In the third approach, one studies ``features'' of persistence diagrams. We will provide a bit more detail about the second and third approaches.

In the second approach, one considers an appropriately defined ``space of persistence diagrams,'' defines a distance function on it, studies geometric properties of this space, and does standard statistical calculations (means, medians, statistical tests, and so on). 
Recall that a persistence diagram (see Fig.~\ref{figure: PH} for an example) is a multiset of points in  $\overline{\mathbb{R}}^2$ and that it gives the same information as a barcode. We now give the following precise definition of a persistence diagram.
\begin{definition}
A \emph{persistence diagram} is a  multiset that  is  the union of a finite multiset of points  in $\overline{\mathbb{R}}^2$  with the multiset of points on the diagonal  $\Delta=\{(x,y)\in \mathbb{R}^2\mid x=y\}$, where each point on the diagonal has infinite multiplicity.
\end{definition}

In this definition, we include all of the points on the diagonal in $\mathbb{R}^2$ with infinite multiplicity. We include the diagonal points for technical reasons; roughly, it is desirable to be able to compare persistence diagrams by studying bijections between their points, and persistence diagrams must thus be sets with the same cardinality.

Given two persistence diagrams $X$ and $Y$, we consider the following general definition of distance between $X$ and $Y$.
\begin{definition}
Let $p \in [1,\infty]$. The {\it $p$th Wasserstein distance} between $X$ and $Y$ is defined as 
\begin{equation*}
	W_p[\dist](X,Y):= \inf_{\phi\colon X\to Y}\left [ 
\sum_{x\in X}\dist[x,\phi(x
)]^p
\right]^{1/p}
\end{equation*}
for $p \in [1,\infty)$
and as 
\begin{equation*}
	W_\infty[\dist](X,Y):= \inf_{\phi\colon X\to Y}\sup_{x\in X}
\dist[x,\phi(x)]
\end{equation*}
for $p = \infty$, where $\dist$ is a metric on $\mathbb{R}^2$ and $\phi$ ranges over all bijections from $X$ to $Y$.
\end{definition}

Usually, one takes $\dist=L_q$ for $q \in [1,\infty]$. One of the most commonly employed distance functions is the {\it bottleneck distance} $W_\infty[L_\infty]$.

The development of statistical analysis on the space of persistence diagrams is an area of ongoing research, and presently there are few tools that can be used in applications. See \cite{MMH11,TMMH14,MTBMMH15} for research in this direction. Until recently, the library $\dionysus$ \cite{dionysus} was the only library to implement computation of the bottleneck and Wasserstein distances (for $d=L_\infty$); the library $\hera$ \cite{hera} implements a new algorithm \cite{KMN16} for the computation of the bottleneck and Wasserstein distances that significantly outperforms the implementation in $\dionysus$.  
The library $\tdapackage$ \cite{tda} (see \cite{FKLM14} for the accompanying tutorial) implements the computation of confidence sets for persistence diagrams that was developed in \cite{FLRWBS14}, distance functions that are robust to noise and outliers \cite{CFLMRW14}, and many more tools for interpreting barcodes. 

The third approach for the development of statistical tools for PH 
consists of mapping the space of persistence diagrams to spaces (e.g., Banach spaces) that are amenable to statistical analysis and machine-learning techniques. Such methods include persistence landscapes \cite{B15}, using the space of algebraic functions \cite{ACC13}, persistence images \cite{CEHKMNPSZ15}, and kernelization techniques \cite{KHNWB15,RHBK15,BMT17,ZVBNB16}. 
See the papers \cite{KNBNH14,SHP16} for applications of persistence landscapes.
The  package $\landscape$  \cite{landscape}  (see \cite{BD15} for the accompanying tutorial) implements the computation of persistence landscapes, as well as many statistical tools that one can apply to persistence landscapes, such as mean, ANOVA,  hypothesis tests, and many more.


\subsection{Stability}\label{SS:stability}

As we mentioned in Section \ref{intro}, PH is useful for applications because it is stable with respect to small perturbations in the input  data. 

The first stability theorem for PH, proven in \cite{CSEH07}, asserts that, under favorable conditions, step (2) in the pipeline in Fig.~\ref{figure:first pipeline} is $1$-Lipschitz with respect to suitable distance functions on filtered complexes and the bottleneck distance for barcodes (see Section \ref{SS:interpretation}). This result was generalized in the papers \cite{CCSGGO09,BS14,BdSS14}. Stability for PH is an active area of research; for an overview of stability results, their history and recent developments, see \cite[Chapter 3]{Ou15}. 


\section{Excursus: Generalized persistence}\label{S:generalised}

One can use the algorithms that we described in Section \ref{S:computation} to compute PH when one has a sequence of complexes with inclusion maps that are all going in the same direction, as in the following diagram:
\begin{equation*}
	\dots \rightarrow K_{i-1}\rightarrow K_i\rightarrow K_{i+1}\rightarrow \dots \,.
\end{equation*}
An algorithm, called the {\tt zigzag algorithm}, for the computation of PH for inclusion maps that do not all go in the same direction, as, e.g., in the diagram
\begin{equation*}
	\dots \rightarrow K_{i-1} \rightarrow K_i\leftarrow K_{i+1}\rightarrow \dots 
\end{equation*}
\noindent
was introduced in  \cite{CSM09}. In the more general setting in which maps are not inclusions, one can still compute PH using the {\tt simplicial map algorithm} \cite{DFW14}.

One may also wish to vary two or more parameters instead of one. This yields multi-filtered simplicial complexes, as, e.g., in the following diagram: 
\[
	\begin{array}{ccccccccc}
&				& \rotatebox{90}{$\dots$} 			& 			&\rotatebox{90}{$\dots$} 		&			&\rotatebox{90}{$\dots$} 		&			&		\\	
&				&\rotatebox{90}{$\rightarrow$}		& 			&\rotatebox{90}{$\rightarrow$}	&			& \rotatebox{90}{$\rightarrow$}	&			&		\\
\dots 			& \rightarrow 	&  K_{j+1,i-1} 					& \rightarrow 	&K_{j+1,i} 				& \rightarrow 				&  K_{j+1,i+1}				& \rightarrow 	& \dots 	\\
&				&\rotatebox{90}{$\rightarrow$}		& 			&\rotatebox{90}{$\rightarrow$}	&			& \rotatebox{90}{$\rightarrow$}	&			&		\\
\dots 			& \rightarrow 	&  K_{j,i-1} 					& \rightarrow 	&K_{j,i} 					& \rightarrow 				&  K_{j,i+1}				& \rightarrow 	& \dots 	\\
&				&\rotatebox{90}{$\rightarrow$}		& 			&\rotatebox{90}{$\rightarrow$}	&			& \rotatebox{90}{$\rightarrow$}	&			&		\\
\dots 			& \rightarrow 	&  K_{j-1,i-1} 					& \rightarrow 	&K_{j-1,i} 					& \rightarrow 				&  K_{j-1,i+1}				& \rightarrow 	& \dots 	\\
&				&\rotatebox{90}{$\rightarrow$}		& 			&\rotatebox{90}{$\rightarrow$}	&			& \rotatebox{90}{$\rightarrow$}	&			&		\\
&				& \rotatebox{90}{$\dots$} 			& 			&\rotatebox{90}{$\dots$} 		&			&\rotatebox{90}{$\dots$} 		&			&		\\	
	\end{array}\,.
\]

\noindent
In this case, one speaks of {\it multi-parameter persistent homology}.  Unfortunately, the Fundamental Theorem of Persistent Homology is no longer valid if one filters with more than one parameter, and there is no such thing as a ``generalized interval.'' 
The topic of multi-parameter persistence is under active research, and several approaches are being studied to extract topological information from multi-filtered simplicial complexes. See \cite{CZ09,Ou15} for the theory of multi-parameter persistent homology, and see \cite{rivet} (and \cite{LW15} for its companion paper) for upcoming software for the visualization of $2$-parameter persistent homology.

\section{Software}\label{S:software}

There are several publicly-available implementations for the computation of PH. We give an overview of the libraries with accompanying peer-reviewed publication and summarize their properties in Table \ref{table: software}. 

The software package $\javaplex$ \cite{javaplex}, which was developed by the computational topology group at Stanford University, is based on the $\plex$ library \cite{plex}, which to our knowledge is the first piece of software to implement the computation of PH. 
$\perseus$ \cite{perseus} was developed to implement Morse-theoretic reductions \cite{MN13} (see Section \ref{SS:reduction}). $\jholes$ \cite{jholes} is a $\java$ library for computing the weight rank clique filtration for weighted undirected networks \cite{PSDV13}. \dionysus \cite{dionysus} is the first software package to implement the {\tt dual algorithm} \cite{dSMVJ11,dSMVJ11b}. $\phat$ \cite{phat} is a library that implements several algorithms and data structures for the fast computation of barcodes, takes a boundary matrix as input, and is the first software to implement a matrix-reduction algorithm that can be executed in parallel. $\dipha$ \cite{dipha}, a spin-off of $\phat$, implements a distributed computation of the matrix-reduction algorithm. $\gudhi$ \cite{gudhi} implements new data structures for simplicial complexes and the boundary matrix. It also implements the {\tt multi-field algorithm}, which allows simultaneous computation of PH over several fields \cite{BM14}. This library is currently under intense development, and a $\python$ interface was just released in the most recent version of the library (namely, Version 2.0.0, whereas the version that we study in our tests is Version 1.3.1). The library $\ripser$ \cite{ripser}, the most recently developed software of the set that we examine, uses several optimizations and shortcuts to speed up the computation of PH with the VR complex. This library does not have an accompanying peer-reviewed publication. However, because it is currently the best-performing (both in terms of memory usage and in terms of wall-time seconds\footnote{``Wall time'' is the amount of elapsed time perceived by a human.}) library for the computation of PH with the VR complex, we include it in our study. The library
$\simppers$ \cite{simppers} implements the {\tt simplicial map algorithm}. Libraries that implement techniques for the statistical interpretation of barcodes include the $\tdapackage$ \cite{tda} and the $\landscape$ \cite{landscape}. $\rivet$, a package for visualizing $2$-parameter persistent homology, is slated to be released soon \cite{rivet}.
We summarize the properties of the libraries that we mentioned in this paragraph in Table \ref{table: software}, and we discuss the performance for a selection of them in Section \ref{S:results} and in the SI. 
For a list of programs, see \url{https://github.com/n-otter/PH-roadmap}.

 \begin{sidewaystable}
\footnotesize

\caption{Overview of existing software for the computation of PH that have an accompanying peer-reviewed publication (and also $\ripser$ \cite{ripser}, because of its performance). The symbol ``$-$'' signifies that the associated feature is not implemented. 
For each software package, we indicate the following items. (a) The language in which it is implemented. (b) The implemented algorithms for the computation of barcodes from the boundary matrix. (c) The coefficient fields for which PH is computed, where the letter $p$ denotes any prime number in the coefficient field $\mathbb{F}_p$. (d) The type of homology computed. (e) The filtered complexes that are computed, where $\VR$ stands for Vietoris--Rips complex, $\W$ stands for the weak witness complex, $\W_\nu$ stands for parametrized witness complexes, $\WRCF$ stands for the weight rank clique filtration,  $\alpha$ stands for the alpha complex, and \v{C} for the \v{C}ech complex.
$\perseus$, $\dipha$, and $\gudhi$ implement 
 the computation of the lower-star filtration \cite{EMP06} 
of a weighted cubical complex; one inputs data in the form of a $d$-dimensional array; the data is then interpreted as a $d$-dimensional cubical complex, and its lower-star filtration is computed. (See the Tutorial in the SI for more details.) Note that $\dipha$ and $\gudhi$ use the efficient representation of cubical complexes presented in \cite{WCV12}, so the size of the cubical complex that is computed by these libraries is smaller than the size of the resulting complex with $\perseus$.
(f) The filtered complexes that one can give as input. $\javaplex$ supports the input of a filtered CW complex for the computation of cellular homology \cite{H02}; in contrast with simplicial complexes, there do not currently exist algorithms to assign a cell complex to point-cloud data. 
(g) Additional features implemented by the library.
 $\javaplex$ supports the computation of some constructions from homological algebra (see \cite{javaplex} for details), and $\perseus$ implements the computation of PH with the VR for points with different ``birth times'' (see Section \ref{SSS:finite metric spaces}). The library $\dionysus$ implements the computation of vineyards \cite{CSM09} and circle-valued functions \cite{dSMVJ11}. Both $\javaplex$ and $\dionysus$ support the output of representatives of homology classes for the intervals in a barcode. 
 (i) Whether visualization of the output is provided. 
}
 \label{table: software}
 \begin{tabular}{|c|c|c|c|c|c|c|c|c|c|}
 \toprule
 Software & \javaplex & \perseus & \jholes & \dionysus & \phat &\dipha & \gudhi & \simppers & \ripser\\\hline
(a) Language & \java & \cpp & \java & \cpp & \cpp& \cpp & \cpp & \cpp & \cpp \\\hline
\begin{tabular}{@{}c@{}} (b) Algorithms\\ for PH \end{tabular}& 
\begin{tabular}{@{}c@{}}standard,\\ dual, \\ zigzag \end{tabular} & 
\begin{tabular}{@{}c@{}} Morse\\ reductions, \\standard \end{tabular} &
\begin{tabular}{@{}c@{}} standard\\ (uses \\ \javaplex) \end{tabular} &
\begin{tabular}{@{}c@{}} standard, \\ dual, \\  zigzag \end{tabular} &
\begin{tabular}{@{}c@{}} standard,  \\ dual,  \\ twist, \\ chunk, \\ spectral seq. \end{tabular}&
\begin{tabular}{@{}c@{}} twist,\\ dual, \\ distributed \end{tabular}& 
\begin{tabular}{@{}c@{}}dual,  \\ multifield \end{tabular}&
\begin{tabular}{@{}c@{}}simplicial  \\ map\end{tabular} &
\begin{tabular}{@{}c@{}} twist,\\ dual \end{tabular}   \\\hline
(c) Coeff. field &
\begin{tabular}{@{}c@{}} $\mathbb{Q}$\,, \\ $\mathbb{F}_p$\end{tabular}&
$\mathbb{F}_2$&
$\mathbb{F}_2$ &
\begin{tabular}{@{}c@{}}$\mathbb{F}_2$ (standard, zigzag),\\ $\mathbb{F}_p$ (dual)\end{tabular}&
$\mathbb{F}_2$&
$\mathbb{F}_2$&
$\mathbb{F}_p$&
$\mathbb{F}_2$ &
$\mathbb{F}_p$ \\\hline
(d) Homology&
\begin{tabular}{@{}c@{}} simplicial,\\ cellular \end{tabular}&
\begin{tabular}{@{}c@{}} simplicial,\\ cubical \end{tabular}&
simplicial&
simplicial&
\begin{tabular}{@{}c@{}} simplicial,\\ cubical \end{tabular}&
\begin{tabular}{@{}c@{}} simplicial,\\ cubical \end{tabular}&
\begin{tabular}{@{}c@{}} simplicial,\\ cubical \end{tabular}&
simplicial &
simplicial\\\hline
\begin{tabular}{@{}c@{}} (e) Filtrations\\ computed \end{tabular}&  
\begin{tabular}{@{}c@{}}$\VR$, \\$\W$,\\$\W_\nu$\\\end{tabular}&
\begin{tabular}{@{}c@{}}  $\VR$,\\ lower star\\ of cubical\\ complex \end{tabular}&
\begin{tabular}{@{}c@{}}  $\WRCF$\\ \end{tabular}&
\begin{tabular}{@{}c@{}}  $\VR$,\\ $\alpha$, \\ \u{C}\end{tabular}&
\begin{tabular}{@{}c@{}}  $-$\\ \end{tabular}&
\begin{tabular}{@{}c@{}}  $\VR$,\\ lower star\\ of cubical \\ complex\end{tabular}&
\begin{tabular}{@{}c@{}}  $\VR$,\\$\alpha$,\\ $\W$ \\lower star\\ of cubical \\complex \end{tabular}&
\begin{tabular}{@{}c@{}} $-$  \\  \end{tabular} &
$\VR$ \\\hline
\begin{tabular}{@{}c@{}} (f) Filtrations\\ as input \end{tabular}&
\begin{tabular}{@{}c@{}} simplicial \\complex, \\\ zigzag, \\ CW \end{tabular}&
\begin{tabular}{@{}c@{}} simplicial\\ complex, \\cubical \\complex \\ \end{tabular}&
$-$ &
\begin{tabular}{@{}c@{}}  simplicial \\complex,\\zigzag\end{tabular}&
\begin{tabular}{@{}c@{}}  boundary \\matrix \\ of simpl. \\complex\end{tabular}&
\begin{tabular}{@{}c@{}}  boundary \\matrix \\ of simpl.\\ complex\end{tabular}&
$-$ &
\begin{tabular}{@{}c@{}}   map of simpl.\\ complexes\end{tabular}&
$-$ \\\hline
\begin{tabular}{@{}c@{}} (g) Additional \\ features \end{tabular}&
\begin{tabular}{@{}c@{}}  Computes \\ some hom. alg.\ \\ constructions, \\ homology\\ generators \end{tabular}&
\begin{tabular}{@{}c@{}}  weighted\\ points \\ for VR \end{tabular}&
$-$&
\begin{tabular}{@{}c@{}}  vineyards,\\ circle-valued \\functions,\\  homology\\ generators  \end{tabular}&
$-$&
$-$&
$-$&
$-$&
$-$\\\hline
(h) Visualization&
barcodes&
\begin{tabular}{@{}c@{}} persistence \\diagram \end{tabular}&
$-$&
$-$&
$-$&
\begin{tabular}{@{}c@{}} persistence \\ diagram \end{tabular}&
$-$&
\begin{tabular}{@{}c@{}} $-$\end{tabular}&
$-$\\
\bottomrule
\end{tabular}
\end{sidewaystable}


 \subsection{Benchmarking}\label{SS:benchmarking}

We benchmark a subset of the currently available open-source libraries with peer-reviewed publication for the computation of PH. To our knowledge, the published open-source libraries are $\jholes$, $\javaplex$, $\perseus$, $\dionysus$, $\phat$, $\dipha$, $\simppers$, and $\gudhi$. To these, we add the library $\ripser$, which is currently the best-performing library to compute PH with the VR complex. To study the performance of the packages, we restrict our attention to the algorithms that are implemented by the largest number of libraries. These are the VR complex and the standard and dual algorithms for the reduction of the boundary matrix. $\phat$ only takes a boundary matrix as input, so it is not possible to conduct a direct comparison of it with the other implementations. However, the fast data structures and algorithms implemented in $\phat$ are also implemented in its spin-off software $\dipha$, which we include in the benchmarking. 
The software $\jholes$ computes PH using the WRCF for weighted undirected networks, and $\simppers$ takes a map of simplicial complexes as input, so these two libraries cannot be compared directly to the other libraries. In the SI, we report benchmarking of some additional features that are implemented by some of the six libraries (i.e., $\javaplex$, $\perseus$, $\dionysus$, $\dipha$, $\gudhi$, and $\ripser$) that we test. 
Specifically, we report results for the computation of PH with cubical complexes for image data sets and the computation of PH with witness, alpha, and \v{C}ech complexes.

We study the software packages $\javaplex$, $\perseus$, $\dionysus$, $\dipha$, $\gudhi$, and $\ripser$ using both synthetic and real-world data from three different perspectives:
\begin{enumerate}
\item Performance measured in CPU seconds and wall-time (i.e., elapsed time) seconds. 
\item Memory required by the process.
\item Maximum size of simplicial complex allowed by the software.
\end{enumerate}


\subsubsection{Data sets}\label{S:data sets}

In this subsection, we describe the data sets that we use for our benchmarking. We use data sets from a variety of different 
mathematical and scientific areas and applications.
In each case, when possible, we use data sets that have already been studied using PH. Our list of data sets is far from complete; we view this list as an initial step towards building a comprehensive collection of benchmarking data sets for PH. 

Data sets \ref{klein}--\ref{fractal} are synthetic: these arise from topology \ref{klein}, stochastic topology \ref{RVR}, dynamical systems \ref{vicsek }, and from an area at the intersection of network theory and fractal geometry \ref{fractal} and which was first used to study connection patterns of the cerebral cortex (see below for details).
Data sets \ref{HIV}--\ref{netw-sc} are from empirical experiments and measurements: they arise from phylogenetics \ref{HIV}--\ref{H3N2}, genomics \ref{genome}, neuroscience \ref{elegans}, image analysis \ref{dragon}, medical imaging \ref{grey}, political science \ref{US}, and scientometrics \ref{netw-sc}.

In each case, these data sets are of one of the following three types: point clouds, weighted undirected networks, and grey-scale digital images. To obtain a point cloud from a real-world weighted undirected network, we compute shortest paths using the inverse of the nonzero weights on edges as distances between nodes (except for the US Congress networks; see below). For the synthetic networks, the values assigned to edges are interpreted as distances between nodes, and we therefore use these values to compute shortest paths. 
We make all processed versions of the data sets that we use in the benchmarking available at \url{https://github.com/n-otter/PH-roadmap/tree/master/data_sets}. We provide the scripts that we used to produce the synthetic data sets at \url{https://github.com/n-otter/PH-roadmap/tree/master/matlab/synthetic_data_sets_scripts}.

\begin{figure}[H]
\begin{center}
\includegraphics[scale=0.27]{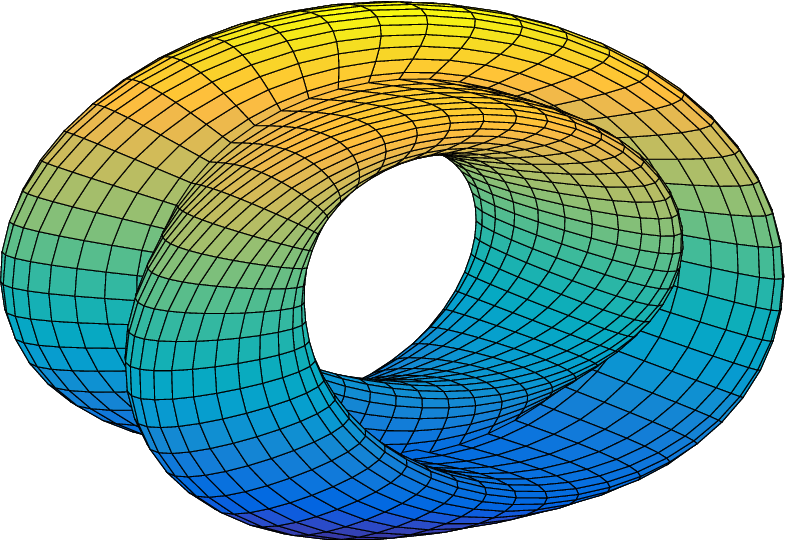} \qquad \qquad \raisebox{0em}{ \includegraphics[scale=0.13]{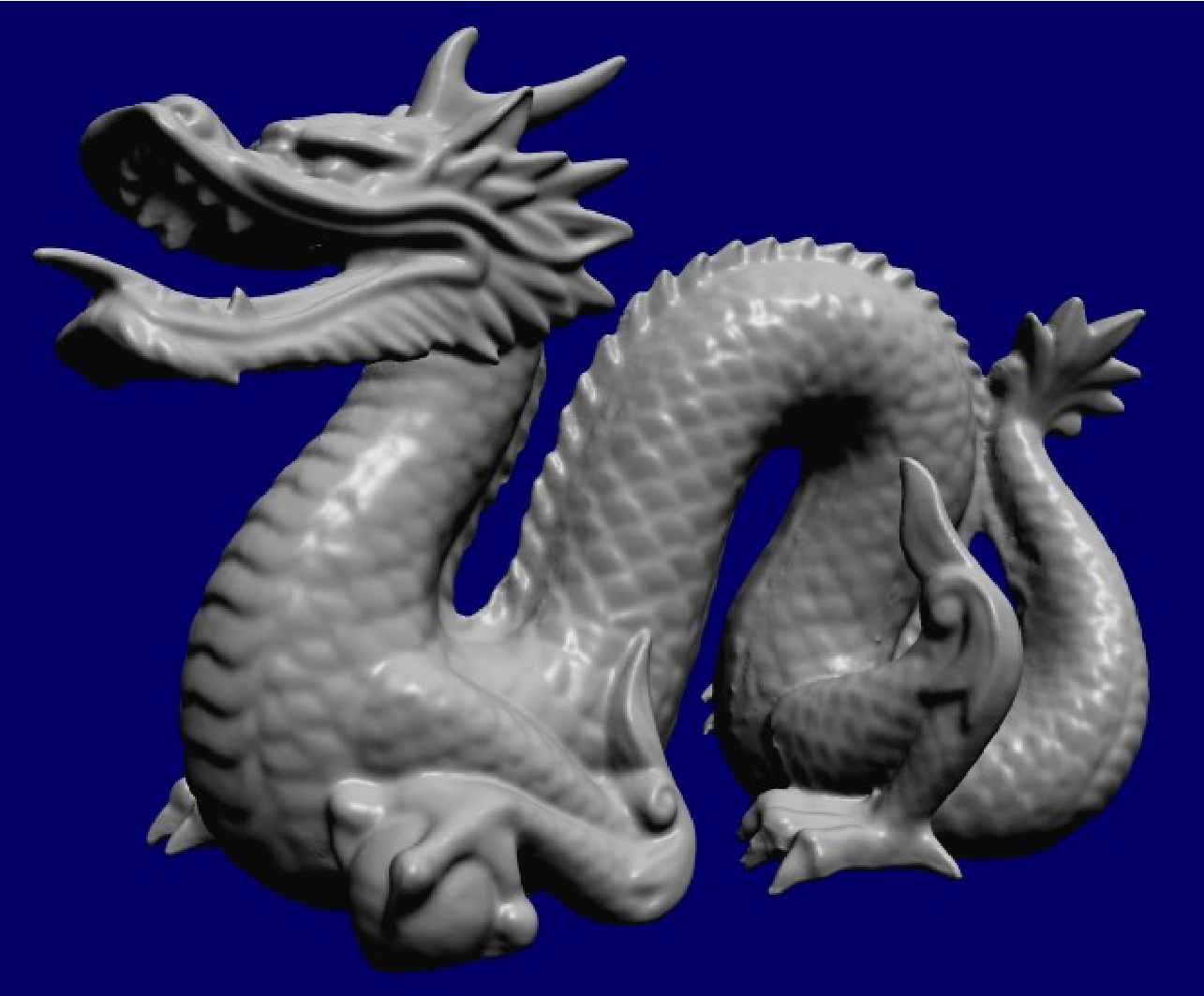} }
\caption{Plot of the image of the figure-$8$ immersion of the Klein bottle and the reconstruction of the Stanford Dragon (retrieved from \cite{Stanford}).
}\label{figure: klein and dragon}
\end{center}
\end{figure}

We now describe all data sets in detail: 
 \begin{enumerate} [label=({\arabic*})]
\item \label{klein} Klein bottle. The Klein bottle is a one-sided nonorientable surface (see Fig.~\ref{figure: klein and dragon}). We linearly sample points from the Klein bottle using its ``figure-8'' immersion in $\mathbb{R}^3$ and size sample of $400$  points. We denote this data set by $\klein$.
Note that the image of the immersion of the Klein bottle does not have the same homotopy type as the original Klein bottle, but they do have the same singular homology\footnote{{\it Singular homology} is a method that assigns to every topological space homology groups encoding invariants of the space, in an analogous way as simplicial homology assigns homology groups to simplicial complexes.  See \cite{H02} for an account of singular homology.} 
 with $\mathbb{F}_2$ coefficients. We have $H_0(B)=\mathbb{F}_2$, $H_1(B)=\mathbb{F}_2\oplus \mathbb{F}_2$, and $H_2(B)=\mathbb{F}_2$, where $B$ denotes the Klein bottle and $H_i(B)$ is the $i$th singular homology group with $\mathbb{F}_2$ coefficients.

\item \label{RVR} 
Random VR complexes (uniform distribution) \cite{K11}. The parameters for this model are positive integers $N$ and $d$; the random VR complex for parameters $N$ and $d$ is the VR complex $\VR_\epsilon(X)$, where $X$ is a set of $N$ points sampled from $\mathbb{R}^d$. (Equivalently, the random VR complex is the clique complex on the random geometric graph $G(N,\epsilon)$ \cite{penrose}.) We sample $N$ points uniformly at random from $[0,1]^d$. We choose $(N,d)=(50,16)$ and we denote this data set by $\random$. The homology of random VR complexes was studied in \cite{K11}.

\item \label{vicsek } Vicsek biological aggregation model. This model was first introduced in \cite{VCBJCS95} and was  studied using PH in \cite{TZH14}. We implement the model in the form in which it appears in \cite{TZH14}. The model describes the motion of a collection of particles that interact in a square with periodic boundary conditions. The parameters for the model are the length $l$ of the side of the square, the initial angle $\theta_0$, the fixed absolute value for the velocity $v_0$, the number of particles $N$, a noise parameter $\eta$, and the number $T$ of time steps. The output of the model is a point cloud in $3$-dimensional Euclidean space in which each point is specified by its position in the $2$-dimensional box and its velocity angle. We run three simulations of the model using the parameter values used in \cite{TZH14}. For each simulation, we choose two point clouds that correspond to two different time frames. 
  See \cite{TZH14} for further details. We denote this data set by $\vicsek $.

\item \label{fractal} Fractal networks. These are self-similar networks introduced in \cite{S06} to investigate whether connection patterns of the cerebral cortex are arranged in self-similar patterns. The parameters for this model are natural numbers $b$, $k$, and $n$. To generate a fractal network, one starts with a fully-connected network on $2^b$ nodes. Two copies of this network are connected to each other so that the ``connection density'' between them is $k^{-1}$, where the connection density is the number of edges between the two copies divided by the number of total possible edges between them. Two copies of the resulting network are connected with connection density $k^{-2}$. One repeats this type of connection process until the network has size $2^n$, but with a decrease in the connnection density by a factor of $1/k$ at each step.

We define distances between nodes in two different ways: (1) uniformly at random, and (2) with linear weight--degree correlations. In the latter, the distance between nodes $i$ and $j$ is distributed as $k_ik_jX$, where $k_i$ is the degree of node $i$ and $X$ is a random variable uniformly distributed on the unit interval. We use the parameters $b=5$, $n=9$, and $k=2$; and we compute PH for the weighted network and for the network in which all adjacent nodes have distance $1$. We denote this data set by $\fract$ and distinguish between the two ways of defining distances between weights using the abbreviations ``r'' for random, and ``l'' for linear. 

 \item \label{HIV} Genomic sequences of the HIV virus. We construct a finite metric space using the independent and concatenated sequences of the three largest genes --- {\tt gag}, {\tt pol}, and {\tt env} --- of the HIV genome. We take 1088 different genomic sequences and compute distances between them by using the Hamming distance. 
We use the aligned sequences studied using PH in \cite{CCR13}. (The authors of that paper retrieved the sequences from \cite{losalamos}.) We denote this data set by $\HIV$.

\item  \label{H3N2} Genomic sequences of H3N2. These are $1000$ different genomic sequences of H3N2 influenza. We compute the Hamming distance between sequences. We use the aligned sequences studied using PH in \cite{CCR13}. We denote this data set by $\influenza$.

\item  \label{dragon} Stanford Dragon graphic. We sample points uniformly at random from 3-dimensional scans of the dragon \cite{Stanford}, whose reconstruction we show in Fig.~\ref{figure: klein and dragon}. The sample sizes contain $1000$ and $2000$ points.  We denote these data sets by $\dragone$ and $\dragtwo$, respectively.

\item \label{elegans} {\it C. elegans} neuronal network. This is a weighted, undirected network in which each node is a neuron and edges represent synapses or gaps junctions. We use the network studied using PH in  \cite{PSDV13}. (The authors of the paper used the data set studied in \cite{WS98}, which first appeared in \cite{WSTB86}.)
Recall that for this example, and also for the other real-world weighted networks (except for the human genome network  and  for the US Congress networks), we convert each nonzero edge weight to a distance by taking its inverse. We denote this data set by $\eleg$.

\item \label{genome} Human genome. A weighted, undirected network representing a sample of the human genome. We use the network studied using PH in \cite{PSDV13}. (The authors  of that paper created the sample using data retrieved from \cite{florida-univ}.) Each node represents a gene, and weighted edges between nodes represent the correlation of the  expression level of the corresponding genes. We define the weight of
an edge as the inverse of the correlation.\footnote{We note that the weight should be the correlation; this issue came to our attention when the paper was in press.}
We denote this data set by $\genome$.

\item \label{grey} Grey-scale image: $3$-dimensional rotational angiography scan of a head with an aneurysm. This data set was used in the benchmarking in \cite{BKRJ14}. This data set is given by a $3$-dimensional array of size $512\times 512\times 512$, with each entry storing an integer that represents the grey value for the corresponding voxel. We retrieved the data set from the repository \cite{volvis}. We denote this data set by $\vertebra$.

\item \label{US} US Congress roll-call voting networks. These two networks (the Senate and House of Representatives from the 104th United States Congress) are constructed using the procedure in \cite{waugh2009} from data compiled by \cite{voteview}. 
In each network, a node is a legislator (Senators in one data set and Representatives in the other), and there is a weighted edge between legislators $i$ and $j$, where the weight $w_{i,j}$ is a number in $[0,1]$ (it is equal to $0$ if and only if legislators $i$ and $j$ never voted the same way on any bill)  given by the number of times the two legislators voted in the same way divided by the total number of bills on which they both voted. See \cite{waugh2009} for further details. We denote the networks from the Senate and House by $\senate$ and $\house$, respectively. The network $\senate$ has $103$ nodes, and the network $\house$ has $445$ nodes. To compute shortest paths, we define the distance between two   nodes $i$ and $j$ to be $1-w_{i,j}$.  
In the 104th Congress, no two politicians voted in the same way on every bill, so we do not have distinct nodes with $0$ distance between them. (This is important, for example, if one wants to apply multidimensional scaling.) 
\item \label{netw-sc} Network of network scientists. This is a weighted undirected network representing the largest connected component of a collaboration network of network scientists \cite{N06}. Nodes represent authors and edges represent collaborations, where weights indicate the number of joint papers. 
The largest connected components consists of $379$ nodes. We denote this data set by $\netwsc$.
 \end{enumerate}

\subsubsection{Machines and compilers}

We tested the libraries on both a cluster and a shared-memory system. The cluster is a Dell Sandybridge cluster, it has $1728$ (i.e., $108 \times 16$) cores of $2.0$GHz Xeon SandyBridge, RAM of $64$ GiB in $80$ nodes and RAM of $128$ GiB in $4$ nodes, and a scratch disk of $20$ TB. It runs the operating system (OS) Red Hat Enterprise Linux $6$. The shared-memory system is an IBM System x3550 M4 server  with $16$ (i.e., $2\times 8$) cores of $3.3$GHz, RAM of $768$ GB, and storage of $3$ TB. It runs the OS Ubuntu 14.04.01.\footnote{Note that we performed the computations for $\gudhi$ and $\ripser$ at a different point in time, during which the shared-memory system was running the OS Ubuntu 16.04.01.} The major difference in running shared algorithms on the shared-memory system versus the distributed-memory system is that each node in the former has much more available RAM than in the latter. (See also the difference in performance between computations on cluster and shared memory system in Tables  \ref{table: performance} and  \ref{table: memory}.) 
To compile $\gudhi$, $\dipha$, $\perseus$ and $\dionysus$ we used the compiler $\gcc$ 4.8.2 on the cluster, and $\gcc$ 4.8.4 on the shared-memory system; for both machines we used the (default) optimization $\textsf{-O3}$. Additionally, we used  $\openmpi$ 1.8.3 for $\dipha$.


\subsubsection{Tests and results}\label{S:results}

We now report the details and results of the tests that we performed. We have made the data sets, header file to measure memory, and other information related to the tests available at \url{https://github.com/n-otter/PH-roadmap}.
Of the six software packages that we study, four implement the computation of the dual algorithm, and four implement the standard algorithm. It is reported in \cite{javaplex} that $\javaplex$ implements the dual algorithm, but the implementation of the algorithm has a bug and gives a wrong output. To our knowledge, this bug has not yet been fixed (at the time of writing),
 and we therefore test only the standard algorithm.

For the computations on the cluster, we compare the libraries running both the dual algorithm and the standard algorithm. The package $\dipha$ is the only one to implement a distributed computation. As a default, we run the software on one node and $16$ cores; we only increase the number of nodes and cores employed when the machine runs out of memory. However, augmenting the number of nodes can make the computations faster (in terms of CPU seconds) for  
complexes of all sizes.\footnote{Based on the results of our tests, we think of small, medium, and large complexes, respectively, as complexes with a size of order of magnitude of up to $10$ million simplices, between $10$ million and $100$ million simplices, and between $100$ million and a billion simplices. 
} 
We see this in our experiments, and it is also discussed in \cite{BKRJ14}. For the other packages, we run the computations on a single node with one core. 
 
For computations on the shared-memory system, we compare the libraries using only the dual algorithm if they implement it, and we otherwise use the standard algorithm. For the shared-memory system, we run all packages (including $\dipha$) on a single core.

In our benchmarking, we report mean computation times and memory measurements. 
In Table \ref{table: performance}, we give the computation times for the different software packages. We measure elapsed and CPU time by using the \texttt{time} function in Linux.  We report computation times with a precision of one second; if a computation took only fractions of a second, we report ``one second'' as the computation time. For space reasons, we report results for a subset of the computations. (In the SI, we tabulate the rest of our computations.) In Table \ref{table: memory}, we report the memory used by the processes in terms of maximum resident set size (RSS); in other words, we give the
maximum amount of real RAM a program has used during its execution. We measure the maximum RSS using the \texttt{getrusage} function in Linux. The header file that we use to measure memory is available at \url{https://github.com/n-otter/PH-roadmap}. In $\dipha$, the measurement of memory is already implemented by the authors of the software. They also use the \texttt{getrusage} function in Linux. 
  The package $\javaplex$ is written in $\java$, and we thus cannot measure its memory as we do for the other packages. However, one can infer memory requirements for this software package using the value of the maximal heap size necessary to perform the computations; we report this value in Table \ref{table: memory}.
In Table \ref{table: failure and max size}, we give the maximum size of the simplicial complex for which we were able to compute PH with each software package in our benchmarkings.

\begin{table}[htbp!]
\footnotesize
\caption{Performance of the software packages measured in wall-time (i.e., elapsed time), and CPU seconds (for the computations running on the cluster). For each data set, we indicate the size of the simplicial complex and the maximum dimension up to which we construct the VR complex.  For all data sets, we construct the filtered VR complex up to the maximum distance between any two points. We indicate the implementation of the standard algorithm using the abbreviation ``st'' following the name of the package, and we indicate the implementation of the dual algorithm using the abbreviation ``d.''  The symbol ``-'' signifies that we were unable to finish computations for this data set, because the machine ran out of memory. $\perseus$ implements only the standard algorithm, and $\gudhi$ and $\ripser$ implement only the dual algorithm. (a,b) We run $\dipha$ on one node and $16$ cores for the data sets $\eleg$, $\klein$, and $\genome$; on $2$ nodes of $16$ cores for the $\HIV$ data set; on $2$ and $3$ nodes of $16$ cores for the dual and standard implementations, respectively, for $\dragtwo$; and on $8$ nodes of $16$ cores for  $\random$. (The maximum number of processes that we could use at any one time was $128$.) (c) We run $\dipha$ on a single core. 
}\label{table: performance}
\begin{center}

 \subfloat[Computations on cluster: wall-time seconds]{

\begin{tabular}{|c|c|c|c|c|c|c|c|c|}
\toprule
Data set  & $\eleg$ & $\klein$& $\HIV$&$\dragtwo$&$\random$ & $\genome$\\\hline
Size of complex& $4.4\times 10^6$&$1.1\times 10^7$ &$2.1\times 10^8$ & $1.3\times 10^9$&$3.1\times 10^9$&$4.5\times 10^8$\\\hline
Max. dim. &2&2&2&2&8&2\\\hline
$\javaplex$ (st)  &84&747&-&-&-&-\\\hline
$\dionysus$ (st)  &474&1830&-&-&-&-\\\hline
$\dipha$ (st) &6&90&1631&142559&-&9110\\\hline
$\perseus$ &543&1978 &-&-&-&-\\\hline
$\dionysus$ (d)  &513&145&-&-&-&-\\\hline
$\dipha$ (d)  &4&6&81&2358&5096&232\\\hline
$\gudhi$ &36& 89&1798&14368&-&4753\\\hline
$\ripser$ & 1&1&2&6&349&3\\
\bottomrule
\end{tabular}
}

 \subfloat[Computations on cluster: CPU seconds]{

\begin{tabular}{|c|c|c|c|c|c|c|c|}
\toprule
Data set  & $\eleg$ & $\klein$& $\HIV$&$\dragtwo$&$\random$ & $\genome$\\\hline
Size of complex& $4.4\times 10^6$&$1.1\times 10^7$ &$2.1\times 10^8$ & $1.3\times 10^9$&$3.1\times 10^9$&$4.5\times 10^8$\\\hline
Max. dim. &2&2&2&2&8&2\\\hline
$\javaplex$ (st) &284&1031&-&-&-&-\\\hline
$\dionysus$ (st)  &473&1824&-&-&-&-\\\hline
$\dipha$ (st) &68&1360&25950&1489615&-&130972\\\hline
$\perseus$ &542 & 1974&-&-&-&-\\\hline
$\dionysus$ (d)  &513&145&-&-&-&-\\\hline
$\dipha$ (d)  &39&73&1276&37572&79691&3622\\\hline
$\gudhi$ &36&88&1794&14351&-&4764\\\hline
$\ripser$ & 1&1&2&5&348&2\\
\bottomrule
\end{tabular}
}

\subfloat[Computations on shared-memory system: wall-time seconds]{
\begin{tabular}{|c|c|c|c|c|c|c|c|c|}
\toprule
Data set  & $\eleg$ & $\klein$& $\HIV$ &$\dragtwo$& $\genome$& $\fractr$ \\\hline
Size of complex&$3.2\times 10^8$&$1.1\times 10^7$&$2.1\times 10^8$&$1.3\times 10^9$&$4.5\times 10^8$&$2.8\times 10^9$\\\hline
Max. dim. &3&2&2&2&2&3\\\hline
$\javaplex$ (st) &13607&1358&43861&-&28064&- \\\hline
$\perseus$ &-&1271&-&-&-&-\\\hline
$\dionysus$ (d) &-&100&142055&35366&-&572764\\\hline
$\dipha$ (d)  &926&13&773&4482&1775&3923\\\hline
$\gudhi$ &381&6&177&1518&442&4590\\\hline
$\ripser$ &2 &1&2&5&3&1517\\
\bottomrule
\end{tabular}
}
\end{center}
\end{table}

\begin{table}[htbp!]
\caption{Memory usage in GB for the computations summarized in Table \ref{table: performance}. For $\javaplex$, we indicate the value of the maximum heap size that was sufficient to perform the computation. The value that we give is an upper bound to the memory usage. For $\dipha$, we indicate the maximum memory used by a single core (considering all cores). See Table \ref{table: performance} for details on the number of cores used. 
}\label{table: memory}

\footnotesize
\begin{center}
 \subfloat[Computations on cluster]{
\begin{tabular}{|c|c|c|c|c|c|c|c|}
\toprule
Data set  & $\eleg$ & $\klein$& $\HIV$&$\dragtwo$&$\random$ & $\genome$\\\hline
Size of complex&$4.4\times10^6$ & $1.1\times 10^7$ & $2.1\times 10^8$&$1.3\times 10^9$&$3.1\times 10^9$&$4.5\times 10^8$\\\hline
Max. dim. &2&2&2&2&8&2\\\hline
$\javaplex$ (st)  &$<5$&$<15$&$>64$&$>64$&$>64$&$>64$\\\hline
$\dionysus$ (st)  &1.3&11.6&-&-&-&-\\\hline
$\dipha$ (st) &0.1&0.2&2.7&4.9&-&4.8\\\hline
$\perseus$ &5.1&12.7&- &-&-&-\\\hline
$\dionysus$ (d)  &0.5&1.1&-&-&-&-\\\hline
$\dipha$ (d)  &0.1&0.2&1.8&13.8&9.6&6.3\\\hline
$\gudhi$ &0.2&0.5& 8.5&62.8&-&21.5\\\hline
$\ripser$ & 0.007&0.02&0.06&0.2&24.7&0.07\\
\bottomrule
\end{tabular}
}

 \subfloat[Computations on shared-memory system]{
\begin{tabular}{|c|c|c|c|c|c|c|c|c|}
\toprule
Data set  & $\eleg$ & $\klein$& $\HIV$&$\dragtwo$& $\genome$& $\fractr$\\\hline
Size of complex&$3.2\times 10^8$&$1.1\times 10^7$&$2.1\times 10^8$&$1.3\times 10^9$&$4.5\times 10^8$&$2.8\times 10^9$\\\hline
Max. dim. &3&2&2&2&2&3\\\hline
$\javaplex$ (st) &$<600$&$<15$&$<700$&$>700$&$<700$&$>700$ \\\hline
$\perseus$ &-&$11.7$&-&-&-&-\\\hline
$\dionysus$ (d) &-&$1.1$&$16.8$&$134.2$&-&268.5\\\hline
$\dipha$ (d)  &$31.2$&$0.9$&$17.7$&$109.5$&$37.3$&276.1\\\hline
$\gudhi$ &15.4&0.5&10.2&62.8&21.4&134.8\\\hline
$\ripser$ &0.2 &0.03&0.07&0.2&0.07&155\\
\bottomrule
\end{tabular}
} 
\end{center}
\end{table}


\begin{table}[htbp!]
\caption{For each software package in (a), we indicate in (b) the maximal size of the simplicial complex supported by it thus far in our tests.
 }\label{table: failure and max size}

\footnotesize
\begin{center}
\begin{tabular}{|c|c|c|c|c|c|c|c|c|}
\toprule
\multirow{2}{*}{(a)} &\javaplex&\multicolumn{2}{|c|}{ $\dionysus$ }&\multicolumn{2}{|c|}{\dipha }&\perseus & \gudhi & \ripser\\ \cline{2-9}
				&st&st&d&st&d&st&d&d\\\hline
(b) &$4.5\cdot 10^8$&$1.1\cdot 10^7$&$2.8\times 10^9$&$1.3\cdot10^{9}$&$3.4\cdot10^{9}$&$1\cdot 10^7$&$3.4\cdot10^{9}$ &$3.4\cdot10^{9}$\\
\bottomrule
\end{tabular}
\end{center}
\end{table}


\subsection{Conclusions from our benchmarking}\label{Conclusion}
 
Our tests suggest that $\ripser$ is the best-performing library currently available for the computation of PH with the Vietoris--Rips complex, and in order of decreasing  performance, that $\gudhi$ and $\dipha$ are the next-best implementations. For the computation of PH with cubical complexes, $\gudhi$ outperforms $\dipha$ by a factor of $3$ to $4$ in terms of memory usage, and $\dipha$ outperforms $\gudhi$ in terms of wall-time seconds by a factor of $1$ to $2$ (when running on one core on a shared-memory system). Both $\dipha$ and $\gudhi$ significantly outperform the implementation in $\perseus$. For the computation of PH with the alpha complex, we did not observe any significant differences in performance between the libraries $\gudhi$ and $\dionysus$. Because the alpha complex has fewer simplices than the other complexes that we used in our tests, further tests with larger data sets may be appropriate in future benchmarkings.

There is a huge disparity between implementations of the dual and standard algorithms. In our benchmarking, the dual implementations outperformed standard ones both in terms of  computation time (with respect to both CPU and wall-time seconds) and in terms of the amount of memory used. This significant difference in performance and memory usage was also revealed for the software package $\dionysus$ in \cite{dSMVJ11b}.

To conclude, in our benchmarking, the fastest software packages were $\ripser$, $\gudhi$, and $\dipha$. 
For small complexes, the software packages $\perseus$ and $\javaplex$ are good choices, because they are the easiest ones to use. (They are the only libraries that need only to be downloaded and are then ``plug-and-play,'' and they have user-friendly interfaces.) Because the library $\javaplex$ implements the computation of a variety of complexes and algorithms, we feel that it is the best software for an initial foray into PH computation.

We now give guidelines for the computation of PH based on our benchmarking. We list several types of data sets in Table \ref{table:data sets} and indicate which software or algorithm that we feel  is best-suited to each one. These guidelines are based on the findings of our benchmarking. Note that one can transform networks into distance matrices, and distance matrices can yield points in Euclidean space using a method such as multi-dimensional scaling. Naturally, given a finite set of points in Euclidean space, we can compute their distance matrix. As we discussed in Section \ref{S:data}, image data can also be considered as a finite metric space, so the indications for distance matrices and points in Euclidean space also apply to image data.

\begin{table}[htbp!]
\footnotesize
\caption{Guidelines for which implementation is best-suited for which data set, based on our benchmarking. Recall that we indicate the implementation of the dual algorithm using the abbreviation ``d'' following the name of a package, and similarly we indicate the implementation of the standard algorithm by ``st''. Note that for smaller data sets one can also use $\javaplex$ to compute PH with VR complexes from points in Euclidean space, and $\perseus$ to compute PH with cubical complexes for image data, and with VR complexes for distance matrices. The library $\jholes$ can only handle  networks with density much less than  $1$. 
}
\label{table:data sets}
\begin{center}
\begin{tabular}{|c|c|c|c|}
\toprule
Data type & Complexes & Suggested software\\\hline
networks & WRCF& $\jholes$\\\hline
image data & cubical&$\gudhi$ or $\dipha$ (st)\\\hline
distance matrix & VR& $\ripser$ \\\hline
distance matrix & W& $\javaplex$\\\hline
points in Euclidean space & VR& $\gudhi$\\\hline
points in Euclidean space &\v{C}  & $\dionysus$\\\hline
points in Euclidean space &$\alpha$ (only in dim $2$ and $3$) & $\dionysus$ ((st) in dim $2$, (d) in dim $3$) or $\gudhi$\\

\bottomrule
\end{tabular}
\end{center}
\end{table}

\section{Future directions}\label{S:considerations}

We conclude by discussing some future directions for the computation of PH. As we saw in Section \ref{S:computation}, much work has been done on step 2 (i.e., going from filtered complexes to barcodes) of the PH pipeline of Fig.~\ref{figure:first pipeline}, and there exist implementations of many fast algorithms for the reduction of the boundary matrix. 
Step 1 (i.e., going from data to a filtered complex) of the PH pipeline 
is an active area of research, but many sparsification techniques (see, e.g., \cite{KS13,S13}) for complexes have yet to be implemented, and more research needs to be done on steps 1 and 3 (i.e., interpreting barcodes; see. e.g., \cite{TMMH14,B15,BK07}) of the PH pipeline. In particular, it is important to develop approaches for statistical analysis of persistent homology.

We believe that there needs to be a community-wide effort to build a library that implements the algorithms and data structures for the computation of PH, and that it should be done in a way that new algorithms and methods can be implemented easily in this framework. This would parallel similar community-wide efforts in fields such as computational algebra and computational geometry,  and libraries such as Macaulay2 \cite{M2}, Sage \cite{sage}, and CGAL \cite{cgal}. 

We also believe that there is a need to create guidelines and benchmark data sets for the test of new algorithms and data structures. The methods and collection of data sets that we used in our benchmarking provide an initial step towards establishing such guidelines and a list of test problems.

\section{Availability of data and materials}
The processed version of the data sets used in the benchmarking and the scripts written for the  tutorial are available at \url{https://github.com/n-otter/PH-roadmap}. The open-source libraries for the computation of PH studied in this paper are available at the references indicated in the associated citations.

\section*{Acknowledgements}
We thank the Rabadan Lab at Columbia University for providing the HIV and H3N2 sequences used in \cite{CCR13} and Giovanni Petri for sharing the data sets used in \cite{PSDV13}. We thank Krishna Bhogaonker, Adrian Clough, Patrizio Frosini, Florian Klimm, Vitaliy Kurlin, Robert MacKay,  James Meiss, Dane Taylor, Leo Speidel, Parker Edwards, and Bernadette Stolz  for helpful comments on a draft of this paper. We also thank the anonymous referees for their many helpful comments. The first author thanks Ulrich Bauer, Michael Lesnick, Hubert Wagner, and Matthew Wright for helpful discussions, and  thanks  Florian Klimm, Vidit Nanda, and Bernadette Stolz for precious advice. 
The authors would like to acknowledge the use of the University of Oxford Advanced Research Computing (ARC) facility (\url{http://dx.doi.org/10.5281/zenodo.22558}) in carrying out some of the computations performed in this work.   The first author thanks the support team at the ARC for their  assistance. NO and PG are grateful for support from the EPSRC grant EP/G065802/1 (The Digital Economy HORIZON Hub). HAH gratefully acknowledges EPSRC Fellowship EP/K041096/1. NO and UT were supported by The Alan Turing Institute through EPSRC grant EP/N510129/1.  NO and HAH were supported by the EPSRC institutional grant D4D01270 BKA1.01.

\bibliographystyle{siam}
\bibliography{roadmap}

\end{document}